\documentclass{amsart}
\pdfoutput=1

\usepackage{amsmath,amsthm,amsfonts,amssymb,amscd}
\usepackage{psfrag,graphicx}
\usepackage{natbib}
\usepackage{tikz}
\usetikzlibrary{trees, fit, positioning}

\usepackage{hyperref}

\theoremstyle{plain}

\newtheorem{Thm}{Theorem}[section]
\newtheorem{Lemma}[Thm]{Lemma}
\newtheorem{Prop}[Thm]{Proposition}

\newtheorem{Def}[Thm] {Definition}
\newtheorem{Rmk}[Thm] {Remark}

\newcommand{\RomanNumeral}[1]{\MakeUppercase{\romannumeral #1}}

\newcommand{\cL}{\mathcal{L}}
\newcommand{\cM}{\mathcal{M}}

\newcommand{\NN}{\mathbb{N}}

\newcommand{\QQ}{\mathbb{Q}}
\newcommand{\RR}{\mathbb{R}}

\newcommand{\ZZ}{\mathbb{Z}}

\begin{document}

\title{On irrationals with Lagrange value exactly $3$}


\date{\today}

\author[Z. Cao]{Zhe Cao}
\address[Zhe Cao]{School of Mathematical Sciences, Nankai University, Tianjin, China.}
\email{zhecao@mail.nankai.edu.cn}

\author[H. Erazo]{Harold Erazo}
\address[Harold Erazo]{IMPA, Estrada Dona Castorina, 110. Rio de Janeiro, Rio de Janeiro-Brazil.}
\email{harold.erazo@impa.br}

\author[C. G. Moreira]{Carlos Gustavo Moreira}
\address[Carlos Gustavo Moreira]{SUSTech International Center for Mathematics, Shenzhen, Guangdong, P. R. China and IMPA, Estrada Dona Castorina, 110. Rio de Janeiro, Rio de Janeiro-Brazil.}
\email{gugu@impa.br}

\thanks{The first author would like to thank the Instituto de Matemática Pura e Aplicada (IMPA) where this work began. The second author was partially supported by CAPES and FAPERJ. The third author was partially supported by CNPq and FAPERJ}

\begin{abstract}
	For $c>0$, let $X_c$ denote the set of $x\in\RR\backslash\QQ$ such that $\left| x-\frac{p}{q} \right|<\frac{1}{cq^2}$ has only finitely many rational solutions $\frac{p}{q}$.
	It is a classical fact, known since the 1950s, that $X_c$ is uncountable for $c>3$ and countable for $c<3$.
	However, the cardinality of $X_3$ does not appear to be present in the literature.
	We prove that $X_3$ is uncountable.
	
	More generally, we show that for any $n\in\NN\cup\{\infty\}$, the set of $x\in\RR\backslash\QQ$ with Lagrange value exactly $3$ and such that $\left| x-\frac{p}{q} \right|<\frac{1}{3q^2}$ has exactly $n$ rational solutions $\frac{p}{q}$ is also uncountable.
\end{abstract}

\maketitle

\section{Introduction}\label{sec1}

\subsection{Lagrange Spectrum and bi-infinite sequences}
The Lagrange spectrum is a closed subset of real line that appears naturally in the study of Diophantine approximation. 
Dirichlet's theorem states that for any $x\in\RR\backslash\QQ$, the inequality $\left| x-\frac{p}{q} \right| < \frac{1}{q^2}$ has infinitely many rational solutions $\frac{p}{q}$. 
Hurwitz improved this result by strengthening the constant in the inequality, so for any $x\in\RR\backslash\QQ$, $\left| x-\frac{p}{q} \right| < \frac{1}{\sqrt{5}q^2}$ also has infinitely many rational solutions $\frac{p}{q}$. 
Furthermore, $\sqrt{5}$ is the best constant that works for all irrational numbers.
However, for some fixed irrational $x$, we can still improve the constant.

For $x\in\RR\backslash\QQ$, the \emph{Lagrange value} (also called the \emph{best constant of Diophantine approximation}) of $x$ is
    $$
    k(x)=\sup\left\{ c>0 \mid \left| x - \frac{p}{q} \right| < \frac{1}{cq^2} \text{has infinitely many rational solutions } \frac{p}{q}\right\}
    $$
where $k$ can be regarded as a map from $\RR\backslash\QQ$ to $(0,\infty]$. The \emph{Lagrange spectrum} is the set of all finite Lagrange values:
    $$
    \cL=\left\{ k(x)<\infty \mid x\in\RR\backslash\QQ \right\}.
    $$
Consider the continued fraction expansion of $x\in\RR\backslash\QQ$:
    $$
    x=[x_0;x_1,x_2\dots]=x_0+\frac{1}{x_1+\frac{1}{x_2+\genfrac{}{}{0pt}{0}{}{\ddots}}}
    $$
and let $\frac{p_n}{q_n}=[x_0;x_1,\dots,x_n]$ be the \emph{$n$-th convergent} of $x$.
We have
    $$
    x-\frac{p_n}{q_n} = \frac{(-1)^{n}}{(\gamma_{n+1}+\eta_{n+1})q_n^2}
    $$
where $\gamma_{n+1}=[x_{n+1};x_{n+2},x_{n+3},\dots]$ , $\eta_{n+1}=[0;x_{n},x_{n-1},\dots,x_1]$. 
A classical result by Legendre states that if $\left| x-\frac{p}{q}\right|<\frac{1}{2q^2}$ , where $\frac{p}{q}\in\QQ$, then $\frac{p}{q}=\frac{p_n}{q_n}$ for some $n\in\NN$. 
By the results above, we have an expression of the Lagrange value of $x$:
    \begin{equation}\label{defk}
       k(x)=\limsup_{n\to\infty}(\gamma_{n+1}+\eta_{n+1}).  
    \end{equation}
In this paper, $\NN=\{0,1,2,\dots\}$ represents the set of non-negative integer numbers and $\NN^*=\{1,2,\dots\}$ represents the set of positive integer numbers.
Notice that through the continued fraction, we can identify $x$ with an infinite sequence $(x_i)_{i\in\NN}$, and we are defining Lagrange value for infinite sequences in (\ref{defk}).
We can extend the definition of Lagrange value to bi-infinite sequences: Consider the symbolic dynamical system $(\Sigma,\sigma)$, where $\Sigma=(\NN^{*})^{\ZZ}$ is the symbolic space (also the space of bi-infinite sequences) and $\sigma: \Sigma\rightarrow\Sigma, (x_i)_{i\in\ZZ}\mapsto (x_{i+1})_{i\in\ZZ}$ is the left shift. 
For any $\underline{x}=(x_i)_{i\in\ZZ}\in\Sigma$ , define the \emph{height function} $\lambda$ on $\Sigma$ as
    $$
    \lambda(\underline{x})=[x_0;x_1,x_2,\dots]+[0;x_{-1},x_{-2},\dots]
    $$
and define the \emph{Lagrange value} of $\underline{x}$ as
    $$
    l(\underline{x})=\limsup_{n\to\infty}{\lambda(\sigma^{n}(\underline{x}))}.
    $$
In this case, we can also define the Lagrange spectrum for bi-infinite sequences, but it coincides with the Lagrange spectrum of irrationals, which is
    $$
    \cL=\left\{l(\underline{x})<\infty \mid \underline{x}\in\Sigma \right\}.
    $$
This is because extending the continued fraction of $x$ to the left will not change the value of the limsup. For example, we can extend $(x_i)_{i\in\NN}$ to a bi-infinite sequence $\underline{x}=(x_i)_{i\in\ZZ}$ by setting $x_{-i}=1$ for $i\in\NN^*$ and we have $k(x)=l(\underline{x})$. 
This is known as the dynamical system description of Lagrange spectrum, attributed to Perron. 

Similarly, we can define the \emph{Markov value} of a bi-infinite sequence $\underline{x}\in\Sigma$ as
	\begin{equation}\label{defm}
    	m(\underline{x})=\sup_{n\in\ZZ}{\lambda(\sigma^{n}(\underline{x})})
	\end{equation}
and define the \emph{Markov spectrum} as
    $$
    \cM=\left\{ m(\underline{x})<\infty \mid \underline{x}\in\Sigma \right\}.
    $$
Both $\cL$ and $\cM$ have complicated geometric properties, and it is interesting to figure out the structure of them. 
We refer readers to \cite{cusick1989markoff} and \cite{booklagrange} for more details about these two spectrums.

\subsection{Statement of the main result}\label{sec1.2}

A classical result by Markov \cite{Markoff1879},\cite{Markoff1880} from 1880's states that $\cL$ and $\cM$ coincide before $3$ and this initial part is exactly a discrete sequence accumulating at $3$, which is
	$$
	\cL\cap(-\infty,3)=\cM\cap(-\infty,3)=\left\{ \sqrt{5},\sqrt{8},\frac{\sqrt{221}}{5},\dots \right\}.
	$$
We will recall the Markov's theorem more explicitly in Subsection \ref{subsec2.1}. 
By Markov's theorem we know that the set $\cL\cap(-\infty,3)$ is countable. 
Moreover, the preimage of an element in $\cL\cap(-\infty,3)$ is an equivalence class containing countable elements (see Theorem \ref{thmmarkov}), so the set $k^{-1}(\cL\cap(-\infty,3))$ is also countable.

On the other hand, it is known since the 1950's that the set $\cL\cap(-\infty,3+\varepsilon)$ is uncountable for any $\varepsilon>0$. The earliest reference to this fact is the one given in the survey \cite{Malyshev1981MarkovAL}, page 775, where it is stated that: “P. G. Kogoniya undertook a systematic study of the spectra to the right of the first accumulation point. He proved that 3 is a condensation point of the Markov spectrum.”

For $c>0$, let
    $$
    X_c=\left\{ x\in\RR\backslash\QQ \mid \left| x-\frac{p}{q} \right| < \frac{1}{cq^2} \text{ has only finitely many rational solutions } \frac{p}{q}\right\}.
    $$
It is well known that $X_c$ is uncountable for $c>3$, and is countable for $c<3$. This is because of the containment
	\begin{equation}\label{containment}
		k^{-1}(\cL\cap(-\infty,c))\subseteq X_c\subseteq k^{-1}(\cL\cap(-\infty,c]).
	\end{equation}
These facts follows from very classical results (\cite{Malyshev1981MarkovAL} and \cite{Markoff1879},\cite{Markoff1880}), but the case $c=3$ is apparently not present in the literature.
In fact, Hurwitz \cite{Hurwitz} claimed that $X_3$ is infinite, and more recently, in \cite[Chapter 8]{gugu} it was claimed that $X_3$ is uncountable, however, both references presented these claims without providing proofs. This paper is filling that gap: we prove that the set $X_3$ is also uncountable.

On the other hand, a related question is to ask which numbers have Lagrange value exactly equal to 3, i.e., the structure of the set $k^{-1}(3)$.
It is well known that $k^{-1}(3)$ is uncountable, for example, consider the numbers:
	\begin{equation}\label{degenerate construction}
		y=[\underbrace{1;1,\dots1}_{l_1},2,2,\underbrace{1,1,\dots1}_{l_2},2,2,\underbrace{1,1,\dots1}_{l_3},2,2,\dots]
	\end{equation}
where $l_1<l_2<l_3<\dots$ is an increasing sequence of positive integers and we can show that those $y$ have Lagrange value $3$.

\begin{Rmk}
	The set $k^{-1}(3)$ has Hausdorff dimension $0$, as a consequence of {\rm\cite[Theorem  1]{Moreira2016GeometricPO}}.
	But the set $k^{-1}(3)$ is dense in $\RR$ because of the fact that Lagrange value does not depend on the initial part of the continued fraction.
\end{Rmk}

More generally, for $n\in\NN\cup\{\infty\}$ we consider the set
	$$
	X_3(n)=k^{-1}(3)\cap\left\{ x\in\RR\backslash\QQ \mid \left| x-\frac{p}{q} \right|<\frac{1}{3q^2} \text{ has exactly } n \text{ rational solutions } \frac{p}{q} \right\},
	$$    
so we have $X_3\cap k^{-1}(3)=\cup_{n\in\NN}X_3(n)$.
Notice that by the containment (\ref{containment}) and the fact that $k^{-1}(\cL\cap(0,3))$ is countable, $X_3$ is uncountable if and only if $X_3\cap k^{-1}(3)$ is uncountable. 
Our first result is

\begin{Thm}\label{thm1.2}
	For each $n\in\NN\cup\{\infty\}$ the set $X_3(n)$ is uncountable.
\end{Thm}

\begin{Rmk}
	Actually all the examples in $(\ref{degenerate construction})$ belong to $X_3(\infty)$, so the main difficulty is to find uncountably many elements in $X_3(n)$ where $n\in\NN$.
\end{Rmk}

\begin{Rmk}	
	An irrational number $x$ is called \emph{attainable} if $ \left| x-\frac{p}{q} \right| < \frac{1}{k(x)q^2} $ has infinitely many rational solutions $\frac{p}{q}$.
	
	The set $X_3(\infty)$ corresponds to the attainable elements in $k^{-1}(3)$.
	More precisely, $x\in X_3(\infty)$ means that $\lvert x-\frac{p}{q}\rvert < \frac{1}{3q^2}$ has infinitely many rational solutions but $\lvert x-\frac{p}{q}\rvert < \frac{1}{(3+\varepsilon)q^2}$ has only finitely many rational solutions for any $\varepsilon>0$.
	
	On the other hand, $x\in X_3(0)$ means that $\lvert x-\frac{p}{q}\rvert \geq \frac{1}{3q^2}$ for all $\frac{p}{q}\in\QQ$ but $\lvert x-\frac{p}{q}\rvert < \frac{1}{(3-\varepsilon)q^2}$ has infinitely many rational solutions for all small $\varepsilon>0$.
\end{Rmk}

The challenge to find $x\in X_3(n)$ for $n\in\NN$ is that the initial portion of the continued fraction of $x$ affects the number of solutions.
For this reason, it is reasonable to use projections of bi-infinite sequences in $m^{-1}(3)$, which have extra information because of the supremum in (\ref{defm}).
In \cite{Bom}, Bombieri gave a characterization of $m^{-1}(3)$ through \emph{renormalization operators}, which we will discuss more details in Section \ref{sec2}.

It is relatively easy to give concrete examples of sequences in $m^{-1}(3)$, but these “easy” examples are periodic at both sides and have Lagrange value less than $3$.
Sequences in $m^{-1}(3)\cap l^{-1}(3)$ are not so direct to describe, nevertheless with the help of Bombieri's characterization we have a criteria (see Lemma \ref{lem4l=3}) to produce quite a lot such elements.

\begin{Rmk}
	In this paper, the set $m^{-1}(3)\cap l^{-1}(3)$ was characterized by renormalizations.
	Interestingly, C. Reutenauer {\rm\cite{Reutenauer2006}} characterized the set $m^{-1}(3)\cap l^{-1}(3)$ through Sturmian sequences, which coincides with the union of the classes $(M_2)$ and $(M_3)$ $($notations in {\rm\cite{Reutenauer2006}}$)$.
	However, we will not use that characterization in this paper.
	
\end{Rmk}

Let
	$$
	\pi: \Sigma=(\NN^{*})^{\ZZ} \rightarrow \Sigma^{+}=(\NN^{*})^{\NN^{*}},\quad \underline{x}=(x_i)_{i\in\ZZ}\mapsto \pi(\underline{x})=(x_i)_{i\in\NN^{*}}
	$$
be the projection of a bi-infinite word to the right part. 
For $n\in\NN\cup\{\infty\}$, let
	$$
	Y_3(n)=\left\{ x\in X_3(n) \mid x=[0; \pi(\underline{x})] \text{ for some } \underline{x}\in m^{-1}(3) \right\}.
	$$
Our main result is
\begin{Thm}\label{mainthm}
	For each $n\in\NN\cup\{\infty\}$ the set $Y_3(n)$ is uncountable.
\end{Thm}

A detailed sketch of proof is provided at the beginning of Section $3$.
We first construct uncountably many elements in the set $Y_3(\infty)$ and $Y_3(0)$.
Then by gluing some carefully chosen finite words at the beginning of elements in $Y_3(0)$, we can map $Y_3(0)$ injectively into $Y_3(n)$ for $0<n<\infty$.

The article is organized as follows: 
In Section \ref{sec2} we define renormalization and establish some properties that are required for proving the main result.
We also give some examples of elements in $m^{-1}(3)\cap l^{-1}(3)$ at the end of Section \ref{sec2}.
In Section \ref{sec3} we give a proof of our main result.
Finally, in Section \ref{sec4} we leave some questions about the structure of a new spectrum $\widetilde{\cM}$, which is related to the classical Lagrange and Markov spectra.

\section{Preliminaries}\label{sec2}

\subsection{Markov's Theorem}\label{subsec2.1}

\begin{Def}
	Let $x,y$ be two irrationals, $x$ is \emph{$GL_2(\ZZ)$-equivalent} to $y$ if there exist integers $a,b,c,d$ with $ad-bc=\pm 1$ such that $x=\frac{a+by}{c+dy}$.
\end{Def}

Here is an equivalent description of $GL_2(\ZZ)$-equivalent by continued fractions (see \cite[Theorem A.1]{Bom} for a proof).

\begin{Prop}
	Irrationals $x,y$ are $GL_2(\ZZ)$-equivalent if and only if their continued fractions eventually coincide.
\end{Prop}

\begin{Def}
	A Markov number $m$ is the largest coordinate of a solution of the Markov equation $x^2+y^2+z^2 = 3xyz$ in the positive integers. The multiplicity of a Markov number $m$ is the number of distinct triples $(x, y, z)$ such that $x, y, z$ are positive integers, $x^2+y^2+z^2 = 3xyz$ and $x <= y <= z = m$.
\end{Def} 

A well known uniqueness conjecture by Frobenius states that all the Markov numbers have multiplicity of  $1$ and the conjecture still remains open.

Markov's theorem gives a complete description of the spectrum $\cL$ and $\cM$ at the part before $3$. Here we state a version in \cite{Bom}.

\begin{Thm}\label{thmmarkov}{\rm$($Markov$)$} 
	$$
	\cL\cap(-\infty,3)=\cM\cap(-\infty,3)=\left\{ \sqrt{5},\sqrt{8},\frac{\sqrt{221}}{5},\dots,\frac{\sqrt{9m^2-4}}{m},\dots \right\}
	$$
	where $m$ is a Markov number.
	
	Moreover, let $(m_r)_{r\geq1}$ be the sequence of the ordered Markov numbers with multiplicity, i.e., a non-decreasing sequence whose terms are all Markov numbers and such that the number of times each term appears is equal to its multiplicity.
	There is a sequence of $GL_2(\ZZ)$-inequivalent irrationals $\gamma_{r}=\frac{a_r+\sqrt{9{m_r}^2-4}}{b_r}$ $(a_r,b_r\in\ZZ)$ such that $k(\gamma_{r})=\frac{\sqrt{9m_r^2-4}}{m_r}$, and every $x\in\RR\setminus\QQ$ with $k(x)<3$ is $GL_2(\ZZ)$-equivalent to some $\gamma_{r}$.
\end{Thm}

\begin{Rmk}
	The continued fractions of quadratics are eventually periodic $($see {\rm\cite[Theorem A.1]{Bom}} for a proof$)$.
	In Markov's Theorem, the period of $\gamma_r$ is given by the Cohn tree $($see Figure {\rm\ref{cohntree}}$)$, i.e., every period of $\gamma_r$ is some vertex of the Cohn tree and every vertex of the Cohn tree is some period of $\gamma_r$.
	Without loss of generality, we can choose the $\gamma_r$ in the theorem to be purely periodic.
\end{Rmk}

\subsection{Structure of Admissible Bi-infinite Words}\label{subsec2.2}

Infinite or bi-infinite sequences can also be regarded as infinite or bi-infinite words (of the same alphabet). 
A word can be finite, infinite, or bi-infinite, so we use underlined letters to represent bi-infinite words, and use the usual letters to represent finite or infinite words.
When referring to infinite words we always mean infinite words to the right.

A \emph{factor} $\omega_1$ of some word $\omega$ is a subword of some consecutive letters inside, denoted as $\omega_1\sqsubseteq\omega$. If $\omega=\omega_1\omega_2$, then $\omega_1$ is a \emph{prefix} of $\omega$, and $\omega_2$ is a \emph{suffix}.

\begin{Def}
    A bi-infinite word $\underline{\omega}\in\Sigma$ is called \emph{admissible} if $m(\underline{\omega})\leq 3$, and \emph{strongly admissible} if $m(\underline{\omega})<3$. 
\end{Def}

If there is one letter bigger than $2$ appearing in a bi-infinite word, then the Markov value will be strictly bigger than $3$, so the only possible letters that could appear in an admissible bi-infinite word are $1$ and $2$.
Moreover, it was shown in \cite[Lemma 9]{Bom} that the number of consecutive $1$s or $2$s in an admissible bi-infinite word is always even or infinite, so an admissible bi-infinite word can be written as a bi-infinite word of alphabet $\{a=22,b=11\}$. 
In \cite[Lemma 11]{Bom}, admissible bi-infinite words were classified into four possible types:
    \begin{itemize}
        \item Constant: $a^{\infty}$, $b^{\infty}$.
        \item Degenerate: $a^{\infty}ba^{\infty}$, $b^{\infty}ab^{\infty}$.
        \item Type \RomanNumeral{1}: $\dots ab^{e_{i-1}}ab^{e_i}ab^{e_{i+1}} \dots$ with every $e_i\geq 1$.
        \item Type \RomanNumeral{2}: $\dots a^{f_{i-1}}ba^{f_i}ba^{f_{i+1}}b \dots$ with every $f_i\geq 1$.
    \end{itemize}
The exponent sequences $(e_i)_{i\in\NN},(f_i)_{i\in\NN}$ are called \emph{characteristic sequences}. 
They satisfy the following combinatorial property (\cite[Lemma 13]{Bom}): for for all $i\in\ZZ$, we have
	\begin{align*}
		(e_i-1,e_{i+1},\dots)\preccurlyeq (e_{i-1},e_{i-2},\dots) \\
		(e_i-1,e_{i-1},\dots)\preccurlyeq (e_{i+1},e_{i+2},\dots)
	\end{align*}
where $\preccurlyeq$ is the lexicographic order.

In Bombieri's paper, he introduced the following renormalization algorithm to characterize admissible words: Let $U,V$ be the Nielsen substitutions of free group $F_2\langle a,b \rangle$, which are
	\begin{align*}
		U : 
		\begin{aligned}
			&a\mapsto ab \\
			&b\mapsto b  \\
		\end{aligned}
		,\quad
		V : 
		\begin{aligned}
			&a\mapsto a \\
			&b\mapsto ab  \\
		\end{aligned}
	\end{align*}
and let $u,v$ be inverse of $U,V$, which are
		\begin{align*}
		u : 
		\begin{aligned}
			&a\mapsto ab^{-1} \\
			&b\mapsto b  \\
		\end{aligned}
		,\quad
		v : 
		\begin{aligned}
			&a\mapsto a \\
			&b\mapsto a^{-1}b . \\
		\end{aligned}
	\end{align*}
We only consider positive words, that is, words with positive exponents. 
If the word is no longer positive after applying $u$ or $v$, then we say that the substitution $u$ or $v$ is not well defined.
A remark is that if $u$ or $v$ is well defined for some bi-infinite word $\underline{\omega}$, then all the exponents of $a$ or $b$ in the word $\underline{\omega}$ should be $1$.

We extend some notations for words of alphabet $\{1,2\}$ to words of $\{a,b\}$.
Let $\chi: a\mapsto 22, b\mapsto 11$ be the substitution.
A bi-infinite word $\underline{\omega}$ of alphabet $\{a,b\}$ is \emph{admissible} if the bi-infinite word $\chi(\underline{\omega})$ is admissible.
We also write continued fractions of alphabet $\{a,b\}$,
for example, when we write the continued fraction $[a,b,\dots]$, essentially we mean the continued fraction $[2;2,1,1,\dots]$.
When referring to length of a word (of any type), we will always mean the length in the alphabet $\{1,2\}$.

\begin{Lemma}\label{cfrac}
	Let $[x_0;x_1,\dots]$ and $[y_0;y_1,\dots]$ be two infinite continued fractions, then $[x_0;x_1,\dots]<[y_0;y_1,\dots]$ if and only if $(-1)^{k}x_k<(-1)^{k}y_k$, where $k$ is the first index such that $x_k\neq y_k$.
\end{Lemma}

We use a vertical bar between two letters to represent a \emph{cut} of word (of any type). 
For example, if $x=x_1x_2x_3\dots$ is an infinite word, the cut at the $i+1$-th position is $x_1\dots x_{i-1}x_i\vert x_{i+1}\dots$. 
We define the \emph{value} of a cut of infinite or bi-infinite words by
\begin{align*}
	\lambda(x_1\dots x_{i-1}x_i\vert x_{i+1}\dots)&=[0;x_i,x_{i-1},\dots,x_1]+[x_{i+1};x_{i+2},\dots] \\
	\lambda(\dots x_{i-1}x_i\vert x_{i+1}\dots)&=[0;x_i,x_{i-1},\dots]+[x_{i+1};x_{i+2},\dots]
\end{align*}
where we use the same notation as the height function $\lambda$, which is an abuse of notation. 
From this point of view, height function is the value of the cut at $0$-th position.

We can also define the notation of transpose of words.
Let $x=x_1x_2\dots x_n$ be a finite word, we define the transpose of $x$ to be $x^T=x_n\dots x_2x_1$.
Similarly, let $x=x_1x_2\dots$ be an infinite word and $\underline{x}=\dots x_{-1}x_0x_{1}\dots$ be a bi-infinite word, their transpose are defined by $x^T=\dots x_2x_1$ (infinite to the left) and $\underline{x}^T=\dots x_{1}x_0x_{-1}\dots$ (transpose at $0$-th position).

\begin{Lemma}\label{lemtranspose}
	Let $\underline{\omega}$ be a bi-infinite word of alphabet $\{1,2\}$ and let $E^T \vert x F$ be a cut of $\underline{\omega}$, where $x\in\{1,2\}$ and $E,F$ are two infinite words of alphabet $\{1,2\}$.
	We have $\lambda(E^T \vert x F)=\lambda(F^T\vert x E)$.
\end{Lemma}

\begin{proof}
	This follows from the equality $[0;E]+[x;F]=[x;E]+[0;F]$.
\end{proof}

We will focus on cuts with value bigger than $3$ or not bigger than $3$.

\begin{Def}\label{defbadcut}
	A cut of an infinite or bi-infinite word is \emph{bad} if the cut has value strictly bigger than $3$, otherwise the cut is \emph{good}.
	
	We can extend the definition to finite words. 
	The cut of finite word is \emph{bad (resp. good)} if any extension of the cut in the alphabet $\{a,b\}$ is still a bad cut $($resp. good cut$)$.
	Otherwise we say that the cut is \emph{indeterminate}.
\end{Def}

In the sequel, we will give examples of cuts in words (of any type). In \cite[Lemma 7]{Bom}, it was stated that the cuts $E^T\vert 1F, E^T2\vert22F$ are good and the cuts $E^T1\vert21F, E^T2\vert212F$ are bad, where $E,F$ are any finite or infinite words.
Actually this is equivalent to say that the cuts $\vert 1, 2\vert 22$ are good and the cuts $1\vert21, 2\vert212$ are bad. Since cuts of the form $E^T\vert 1F$ are always good because $\lambda(E^T\vert 1F)=[1;F]+[0;E]<3$, we adopt the following abuse of notation:
\begin{equation*}
    \lambda(E^Ta|b F) = [2;2,E]+[0;1,1,F] = \lambda(E^T2|2b F).
\end{equation*}
Since all the words appearing in this paper can be written in the alphabet $\{a,b\}$, this notation is particularly convenient for applying the following lexicographic comparison criteria.

Same with the lexicographic order of words of alphabet $\{1,2\}$, we can also define a natural lexicographic order on words of alphabet $\{a,b\}$.
The following lemma is from \cite[Lemma 8,10]{Bom}.

\begin{Lemma}\label{critcut}
	A bi-infinite word $\underline{\omega}$ of alphabet $\{1,2\}$ is admissible if and only if 
	\begin{enumerate}
		\item $\underline{\omega}$ can be written as an infinite word of alphabet $\{a,b\}$.
		\item Every cut $E^Tb\vert aF$ of either $\underline{\omega}$ or $\underline{\omega}^T$ satisfies $E\preccurlyeq F$, where $E,F$ are two infinite words of alphabet $\{a,b\}$.
	\end{enumerate}
\end{Lemma}

\begin{proof}
	This follows from Lemma \ref{cfrac}, Lemma \ref{lemtranspose} and the equality $[2;2,\omega]+[0;1,1,\omega]=3$, where $\omega$ is a word of alphabet $\{1,2\}$ (finite or infinite).
\end{proof}

As consequence we can characterize bad and good cuts in bi-infinite words.

\begin{Prop}\label{critbadcut}
	Let $\omega$ be a finite word of alphabet $\{a,b\}$. Then the cuts $a\omega^Tb\vert a\omega b $, $b\omega^Ta\vert b\omega a$ are good and the cuts $b\omega^Tb\vert a\omega a$, $a\omega^Ta\vert b\omega b$ are bad.
\end{Prop}

\begin{proof}
	By Lemma \ref{lemtranspose} we know that the cut $b\omega^T2\vert 2 b\omega a$ can be identified with the cut $a\omega^Tb\vert a\omega b $ and the cut $a\omega^T2\vert 2b\omega b$ can be identified with the cut $b\omega^Tb\vert a\omega a$.
	Then the Proposition follows from Lemma \ref{critcut}.
\end{proof}

We give more examples of good cuts and bad cuts. By Proposition \ref{critbadcut}, the cut $bb\vert aab$ is a bad cut and the cut $aab\vert aab$ is a good cut. However, we can't say anything about the cut $b\vert aab$ because it depends on what is to the left, so the cut $b\vert aab$ is indeterminate.

For every $n\in\NN\cup\{\infty\}$, the elements in $X_3(n)$ are infinite words with Lagrange value $3$ and with exactly $n$ bad cuts.
As was stated in Section \ref{sec1}, our main idea is to take projections of elements in $m^{-1}(3)\cap l^{-1}(3)$, and we have to control the number of bad cuts inside the projections.

Good and bad cuts of a finite word $\omega$ will remain good and bad for any extension (finite, infinite or bi-infinite) of $\omega$.
However, the same does not hold for indeterminate cuts.
For example, the cut $a\vert baab$ is indeterminate since $baa\vert baab$ is a good cut but $aaaa\vert baab$ is a bad cut.
If an infinite word begins with $a\vert baab\dots$, then it must be a bad cut for the infinite word because of the inequality
	$$
	\begin{aligned}
		\lambda(a\vert baab\dots)&=[0;b,a,a,b,\dots]+[2;2]\\
		&>[0;b,a,a,b,\dots]+[2;2,a,a,a]=\lambda(aaaa\vert baab\dots)>3
	\end{aligned}
	$$
where we are using the fact that the length before the cut is odd.

An important observation is that good cuts for a bi-infinite word $\underline{\omega}$ can change to bad for the infinite word $\pi(\underline{\omega})$.
For example, the bi-infinite word $(abb)^{\infty}$ has only good cuts (since it has Markov value less than $3$) but the projection $abbabb\dots$ has bad cuts, because
	$$
		\lambda(a\vert bb\dots)=[0;b,b,\dots]+[2;2]>[0;b,b,\dots]+[2;2,a]=\lambda(aa\vert bb\dots)>3.
	$$
This is the main idea behind creating bad cuts in infinite words.

On the other side, odd cuts remain good under projections, and this will help us control the number of bad cuts in an infinite word.
For an infinite word $\omega=\omega_1\omega_2\dots$, \emph{odd cuts} means cuts at odd positions, that is, the cut $\omega_1\omega_2\dots\omega_{i-1}\vert\omega_i\omega_{i+1}\dots$ is an odd cut if $i$ is odd.
Let $\pi(m^{-1}(3))=\{\pi(\underline{\omega})\mid \underline{\omega}\in m^{-1}(3)\}$ be the set of infinite words which are projections of elements in $m^{-1}(3)$.

\begin{Prop}\label{oddgood}
	Elements in $\pi(m^{-1}(3))$ will have odd good cuts.
\end{Prop}

\begin{proof}
	Given $\underline{\omega}=(\omega_i)_{i\in\ZZ}\in m^{-1}(3)$, we prove that all the odd cuts in $\pi(\underline{\omega})=(\omega_i)_{i\in\NN}$ are good.
	This mainly follows from the inequality
	\begin{align*}
		\lambda(\omega_1\omega_2\dots\omega_{i-1}\vert\omega_i\omega_{i+1}\dots)&=[\omega_i;\omega_{i+1},\dots]+[0;\omega_{i-1},\omega_{i-2},\dots,\omega_1]\\
		&<[\omega_i;\omega_{i+1},\dots]+[0;\omega_{i-1},\omega_{i-2},\dots,\omega_1,\omega_0,\omega_{-1},\dots]\\
		&\leq m(\underline{\omega})=3
	\end{align*}
	where $i$ is an odd position.
\end{proof}

\subsection{Renormalization}\label{subsec2.3}

We state an important lemma, which is essentially a restatement of \cite[Lemma 14]{Bom}.

\begin{Lemma}{\rm (Renormalization algorithm)}\label{renorm algorithm}
	Let $\underline{\omega}$ be a non-constant admissible bi-infinite word, then there exist $\sigma\in\{u,v\}$ such that $\underline{\omega}^{\sigma}$ is well defined and still an admissible bi-infinite word. Moreover, the choice of $\sigma$ is unique if $\underline{\omega}\neq\underline{\omega}_0=\dots ababab \dots$.
	This procedure is called {\rm (}one-step{\rm )} renormalization algorithm.
\end{Lemma} 

Define the renormalization algorithm as follows: given a bi-infinite word $\underline{\omega}$, we can keep applying the one-step renormalization algorithm whenever it is well defined, that is, we can keep choosing a sequence (possibly finite) of $\sigma_i\in\{u,v\}$ such that every $\underline{\omega}^{\sigma_1},\underline{\omega}^{\sigma_1\sigma_2},\underline{\omega}^{\sigma_1\sigma_2\sigma_3},\dots $ is well defined. We define the renormalization algorithm to stop whenever it reaches $\underline{\omega}_0$ or the action $u$ or $v$ is not well defined. 

In the proof of \cite[Theorem 15]{Bom}, Bombieri characterized the admissible sequences by renormalization algorithm in the following sense:

\begin{Prop}\label{prop renorm}
	A non-constant bi-infinite word $\underline{\omega}$ is strongly admissible if and only if the renormalization algorithm stops at some finite step and reaches $\underline{\omega}_0$, and a non-constant bi-infinite word has Markov value exactly $3$ if and only if the renormalization algorithm never stops.
\end{Prop}

We also say that $\underline{\omega}$ is \emph{infinitely renormalizable} if the renormalization algorithm never stops, and \emph{finitely renormalizable} if the renormalizable algorithm stops at some finite step and reaches $\underline{\omega}_0$.

Bombieri also gave a nice characterization of strongly admissible bi-infinite words in \cite[Theorem 15]{Bom}.

\begin{Prop}
	The strongly admissible bi-infinite words are periodic and the periods are given by the \emph{Cohn tree} $T$ $($see Figure $\ref{cohntree})$.
\end{Prop} 

The Cohn tree $T$ is an infinite binary tree generated by the following procedure: start with the root $ab$ and keep applying $U$ to the left and applying $V$ to the right.
The set of vertices of the Cohn tree, denoted as $P$, gives a full description of periods of strongly admissible bi-infinite words.

\begin{figure}[h]
	\centering
	\begin{tikzpicture}[
		level 1/.style={sibling distance=6cm},
		level 2/.style={sibling distance=3.5cm},
		level 3/.style={sibling distance=1.5cm},
		edge from parent/.style={draw, -},
		every node/.style={font=\footnotesize} % Adjust font size if necessary
		]
		\node (A) {$ab$}
		child {node (B) {$abb$}
			child {node (D) {$abbb$}
				child {node (H) {$abbbb$}}
				child {node (I) {$aababab$}}
			}
			child {node (E) {$aabab$}
				child {node (J) {$ababbabb$}}
				child {node (K) {$\cdots$}}
			}
		}
		child {node (C) {$aab$}
			child {node (F) {$ababb$}
				child {node (L) {$\cdots$}}
				child {node (M) {$aabaabab$}}
			}
			child {node (G) {$aaab$}
				child {node (N) {$abababb$}}
				child {node (O) {$aaaab$}}
			}
		};
				
		\path (A) -- (B) node[midway, above] {$U$};
		\path (A) -- (C) node[midway, above] {$V$};
		\path (B) -- (D) node[midway, above] {$U$};
		\path (B) -- (E) node[midway, above] {$V$};
		\path (C) -- (F) node[midway, above] {$U$};
		\path (C) -- (G) node[midway, above] {$V$};
		\path (D) -- (H) node[midway, above, xshift=-0.07cm] {$U$};
		\path (D) -- (I) node[midway, above, xshift=0.07cm] {$V$};
		\path (E) -- (J) node[midway, above, xshift=-0.07cm] {$U$};
		\path (E) -- (K) node[midway, above, xshift=0.07cm] {$V$};
		\path (F) -- (L) node[midway, above, xshift=-0.07cm] {$U$};
		\path (F) -- (M) node[midway, above, xshift=0.07cm] {$V$};
		\path (G) -- (N) node[midway, above, xshift=-0.07cm] {$U$};
		\path (G) -- (O) node[midway, above, xshift=0.07cm] {$V$};
		
	\end{tikzpicture}
	\caption{$T$: The Cohn Tree.}
	\label{cohntree}
\end{figure}

However, it is not so immediate to give examples of bi-infinite words with Markov value exactly $3$.
We will give some examples later in Subsection \ref{subsec2.3}.

A natural way to understand the renormalization is to consider whether a bi-infinite word can be written in some alphabet $\{\alpha,\beta\}$.
For example, let $(\alpha,\beta)$ be some word pair of alphabet $\{a,b\}$, we can define actions of $U,V$ on $(\alpha,\beta)$ by substituting by coordinates:
$$
U(\alpha,\beta)=(\alpha^{U},\beta^{U}),\quad V(\alpha,\beta)=(\alpha^{V},\beta^{V}).
$$
The renormalization algorithm (Lemma \ref{renorm algorithm}) can be restated as follows: Starting with $(\alpha_0,\beta_0)=(a,b)$, if a bi-infinite word $\underline{\omega}$ can be written as a non-constant admissible word in alphabet $(\alpha_n,\beta_n)$, then it can also be written as an admissible word in alphabet $(\alpha_{n+1},\beta_{n+1})=U(\alpha_n,\beta_n)$ or $V(\alpha_n,\beta_n)$.

We consider another renormalization which is essentially equivalent to Bombieri's renormalization, as we will see in Proposition \ref{relbtrenorm}.
Given a word pair $(\alpha,\beta)$ of alphabet $\{a,b\}$, let $\overline{U},\overline{V}$ be the following operations: 
	$$
	\overline{U}: (\alpha,\beta)\mapsto (\alpha\beta,\beta),\quad \overline{V}: (\alpha,\beta)\mapsto (\alpha,\alpha\beta).
	$$

\begin{figure}[h]
	\centering
	\begin{tikzpicture}[
		level 1/.style={sibling distance=6cm},
		level 2/.style={sibling distance=3.5cm},
		level 3/.style={sibling distance=1.5cm},
		edge from parent/.style={draw, -},
		every node/.style={font=\footnotesize} % Adjust font size if necessary
		]
		\node (A) {$(a,b)$}
		child {node (B) {$(ab,b)$}
			child {node (D) {$(abb,b)$}
				child {node (H) {$(abbb,b)$}}
				child {node (I) {$(abb,abbb)$}}
			}
			child {node (E) {$(ab,abb)$}
				child {node (J) {$(ababb,abb)$}}
				child {node (K) {$\cdots$}}
			}
		}
		child {node (C) {$(a,ab)$}
			child {node (F) {$(aab,ab)$}
				child {node (L) {$\cdots$}}
				child {node (M) {$(aab,aabab)$}}
			}
			child {node (G) {$(a,aab)$}
				child {node (N) {$(aaab,aab)$}}
				child {node (O) {$(a,aaab)$}}
			}
		};
		
		\path (A) -- (B) node[midway, above] {$\overline{U}$};
		\path (A) -- (C) node[midway, above] {$\overline{V}$};
		\path (B) -- (D) node[midway, above] {$\overline{U}$};
		\path (B) -- (E) node[midway, above] {$\overline{V}$};
		\path (C) -- (F) node[midway, above] {$\overline{U}$};
		\path (C) -- (G) node[midway, above] {$\overline{V}$};
		\path (D) -- (H) node[midway, above, xshift=-0.07cm] {$\overline{U}$};
		\path (D) -- (I) node[midway, above, xshift=0.07cm] {$\overline{V}$};
		\path (E) -- (J) node[midway, above, xshift=-0.07cm] {$\overline{U}$};
		\path (E) -- (K) node[midway, above, xshift=0.07cm] {$\overline{V}$};
		\path (F) -- (L) node[midway, above, xshift=-0.07cm] {$\overline{U}$};
		\path (F) -- (M) node[midway, above, xshift=0.07cm] {$\overline{V}$};
		\path (G) -- (N) node[midway, above, xshift=-0.07cm] {$\overline{U}$};
		\path (G) -- (O) node[midway, above, xshift=0.07cm] {$\overline{V}$};

	\end{tikzpicture}
    \caption{$\overline{T}$: Tree of alphabets induced by $\overline{U}$ and $\overline{V}$.}
\end{figure}

If we start with word pair $(a,b)$ and keep applying $\overline{U}$ and $\overline{V}$ to it, we will get an infinite binary tree, denote as $\overline{T}$.
Let $\overline{P}$ be the set of vertices of $\overline{T}$ and $\overline{P}_n$ be the set of vertices of $\overline{T}$ with distance exactly $n$ to the root $(a,b)$.

With the operators $\overline{U},\overline{V}$, we can restate the renormalization algorithm: 
Start with $(\alpha_0,\beta_0)=(a,b)$, if a bi-infinite word $\underline{\omega}$ can be written as a non-constant admissible word in some alphabet $(\alpha_n,\beta_n)\in \overline{P}_n$, then it can also be written as an admissible word in some alphabet $(\alpha_{n+1},\beta_{n+1})=\overline{U}(\alpha_n,\beta_n)$ or $\overline{V}(\alpha_n,\beta_n)$ in $\overline{P}_{n+1}$. 
Moreover, if $\underline{\omega}\neq \dots\alpha_n\beta_n\alpha_n\beta_n\alpha_n\beta_n\dots$, then the choice is unique: If $\underline{\omega}$ is of type \RomanNumeral{1} or degenerate case $\beta_n^{\infty}\alpha_n\beta_n^{\infty}$ (resp. type \RomanNumeral{2} or degenerate case $\alpha_n^{\infty}\beta_n\alpha_n^{\infty}$) then the next operator will be $\overline{U}$ (resp. $\overline{V}$).
So the Proposition \ref{prop renorm} can be restated as follows:

\begin{Prop}\label{prop renorm 2}
	A bi-infinite word $\underline{\omega}$ is strongly admissible if and only if it can be written as a constant word in some alphabet $(\alpha,\beta)\in \overline{P}$, and $\underline{\omega}$ has Markov value exactly $3$ if and only if for all $n\in\NN$, the word $\underline{\omega}$ can be written as a non-constant word for some alphabet $(\alpha_n,\beta_n)\in\overline{P}_n$.
\end{Prop}

\begin{proof}
	Actually Bombieri \cite[Theorem 15]{Bom} only gave a proof that having Markov value exactly $3$ implies infinitely renormalizable. 
	Here we give a short proof of infinitely renormalizable implies Markov value $3$, assuming Markov's theorem.
	
	If $\underline{\omega}$ has Markov value smaller than $3$, it will be periodic and finitely renormalizable, which contradict to the fact that $\underline{\omega}$ is infinitely renormalizable.
	On the other hand, if $\underline{\omega}$ has Markov value bigger than $3$, then by Lemma \ref{critcut} there exist a bad cut of form $b\omega^T b\vert a\omega a$ in $\underline{\omega}$ or $\underline{\omega}^T$.
	
	If both operators $\overline{U},\overline{V}$ appear infinitely many times in the renormalization of $\underline{\omega}$, the length of both letters $\alpha_n,\beta_n$  will go to infinity, so the subword $b\omega^Tb\vert a\omega a$ will be a factor of $\alpha_N$ or $\beta_N$ for some $N$ big enough.
	This contradicts to the fact that $\alpha_N,\beta_N$ are periods of some strongly admissible words, which do not contain bad cuts.
	
	If one of $\overline{U},\overline{V}$ appears only finitely many times, without loss of generality we assume that $\overline{U}$ only appear at the first $N_1$ steps.
	We claim that under $N_1$ steps of renormalization, the word $\underline{\omega}$ will reach the degenerate case $\alpha_{N_1}^{\infty}\beta_{N_1}\alpha_{N_1}^{\infty}$.
	
	First notice that the word $\underline{\omega}$ can always be written in alphabet $\{\alpha_n,\beta_n\}$ for any $n\in\NN$ and the operators $\overline{U},\overline{V}$ are always well defined, so all the exponents of $\alpha_n$ or $\beta_n$ in $\underline{\omega}$ are $1$.
	In our case after the first $N_1$ steps, the exponents of $\beta_n (n\geq N_1)$ are always $1$.
	We are going to show that $\underline{\omega}=\alpha_{N_1}^{\infty}\beta_{N_1}\alpha_{N_1}^{\infty}$ is the only possible option.
	If there exist some $\alpha_{N_1}$ with finite exponent, then under some finite step $N_2>N_1$ of renormalization the exponent will reduce to $0$. 
	But that will produce some $\beta_{N_2}$ such that the exponent is bigger than $1$, which means that in the step $N_2$ the word will have all exponents $1$ on the letter $\alpha_{N_2}$.
	Moreover, in the $N_2$ step of renormalization the operator will again be $\overline{U}$, which  contradicts the assumption.
	So the exponent of $\alpha_{N_1}$ in $\underline{\omega}$ should be infinity, and the only possible option is the degenerate case $\alpha_{N_1}^{\infty}\beta_{N_1}\alpha_{N_1}^{\infty}$. 
	
	However, this again contradicts with the fact that $\underline{\omega}$ has Markov value bigger than $3$ since in this case the Markov value of $\underline{\omega}$ is going to be $3$.
\end{proof}

We call the renormalization induced by $U,V$ \emph{interior renormalization} and the renormalization induced by $\overline{U},\overline{V}$ \emph{exterior renormalization}.
The operators $U,V,\overline{U},\overline{V}$ are called \emph{renormalization operators}.
The two renormalizations are equivalent in the following sense:

\begin{Prop}\label{relbtrenorm}
	Given $(\alpha_n,\beta_n)\in\overline{P}_n$ a word pair, there exist a finite word $\overline{R}=\overline{R}_1\overline{R}_2\dots\overline{R}_n$ of alphabet $\{\overline{U},\overline{V}\}$ such that $(\alpha_n,\beta_n)=\overline{R}(a,b)$. Let $R_1,R_2,\dots R_n$ be the corresponding interior renormalization operators $($simply remove the bar$)$ and $R=R_1R_2\dots R_n$.
	We have the following equality
	$$
	(\alpha_n,\beta_n)=\overline{R}(a,b)=R^T(a,b).
	$$
\end{Prop}

\begin{proof}
	We use induction to prove the proposition.
	For the case $n=1$ the proposition holds because we have $(ab,b)=\overline{U}(a,b)=U(a,b)$ and $(a,ab)=\overline{V}(a,b)=V(a,b)$.
	
	Now suppose that the proposition holds for the case $n=k$, we prove for the case $n=k+1$. Suppose that $(\alpha_{k+1},\beta_{k+1})=\overline{R}_1\overline{R}_2\dots\overline{R}_{k+1}(a,b)=\overline{R}_1(\alpha_k,\beta_k)\in\overline{P}_{k+1}$. 
	Let $\overline{R}=\overline{R}_2\dots\overline{R}_{k+1}$ and let $R_1,R_2,\dots R_{k+1}, R=R_2\dots R_{k+1}$ be the corresponding interior renormalization operator.
	By induction we have $(\alpha_k,\beta_k)=\overline{R}(a,b)=R^T(a,b)$, so we have
	\begin{align*}
		\overline{U} \overline{R}(a,b)&=(\alpha_k\beta_k,\beta_k)=(a^{R^T}b^{R^T},b^{R^T})=((ab)^{R^T},b^{R^T})=R^T U(a,b)=(UR)^T(a,b) \\
		\overline{V} \overline{R}(a,b)&=(\alpha_k,\alpha_k\beta_k)=(a^{R^T},a^{R^T}b^{R^T})=(a^{R^T},(ab)^{R^T})=R^T V(a,b)=(VR)^T(a,b)
	\end{align*}
	which finishes the induction.
\end{proof}

\subsection{Properties of Renormalization}\label{subsec2.4}

We state a lemma from (\cite[Lemma 3.8]{Har}) and fix some notations.
Let $\omega=\omega_1\omega_2\dots\omega_n$ be a finite word of alphabet $\{a,b\}$. 
Let $\omega^+=\omega_2\dots\omega_n$ be the word obtained by removing the first letter, similarly let $\omega^-=\omega_1\dots\omega_{n-1}$ be the word obtained by removing the last letter.
If $\omega$ starts with $a$ then we define $\omega^b=b\omega^+$ be the word obtained by replacing the first letter $a$ with $b$.
Similarly if $\omega$ ends with $b$, then we define $\omega_a=\omega^-a$ be the word obtained by replacing the last letter $b$ with $a$.

\begin{Lemma}\label{Lem3.8Har}
	For every $(\alpha,\beta)\in \overline{P}$, $\alpha$ starts with $a$ and $\beta$ ends with $b$.
	We always have the equalities $\alpha\beta=\beta_a\alpha^b$, $\alpha^b\beta=\beta^T\alpha^b$ and $\alpha\beta_a=\beta_a\alpha^T$.
	As a consequence, words $\alpha^k\beta,\alpha\beta^k$, with $k\geq 1$, start with $\beta_a$ and end with $\alpha^b$.
\end{Lemma}

\begin{proof}
	One can find the proof in \cite[Lemma 3.8]{Har}, which was proved by induction.
	Here we give another proof based on the equality $\alpha\beta=\beta_a\alpha^b$.
	The equality $\alpha\beta=\beta_a\alpha^b$ can be proved by induction: 
	If $(\alpha^\prime,\beta^\prime)=(\alpha,\alpha\beta)$ then we have $\alpha^\prime\beta^\prime=\alpha(\alpha\beta)=\alpha\beta_a\alpha^b=(\alpha\beta)_a\alpha^b=\beta^\prime_a{\alpha^\prime}^b$.
	The case for $(\alpha^\prime,\beta^\prime)=(\alpha\beta,\beta)$ is analogous.
	
	The lemma holds trivially for the case $(\alpha,\beta)=(a,b)$, so we assume that $(\alpha,\beta)\neq(a,b)$ in the following proof.
	Notice that we have $\alpha^b\beta=(\alpha\beta)^b=(\beta_a\alpha^b)^b=(\beta_a)^b\alpha^b=\beta^T\alpha^b$, where $(\beta_a)^b=\beta^T$ follows from the fact that $\beta$ is of the form $a\theta b$ with $\theta$ palindrome (see Lemma \ref{lem1}).
	So for any $k\geq 1$ we have
	$\alpha\beta^k=\beta_a\alpha^b\beta^{k-1}=\beta_a(\beta^T)^{k-1}\alpha^b$, which starts with $\beta_a$ and ends with $\alpha_b$.
	Analogously we have $\alpha\beta_a=\beta_a\alpha^T$ and $\alpha^k\beta=\beta_a(\alpha^T)^{k-1}\alpha^b$, which also starts with $\beta_a$ and ends with $\alpha^b$.
\end{proof}

The following proposition establishes that the representation of any word in the alphabet $(\alpha,\beta)\in\overline{P}$ is unique.

\begin{Prop}\label{propuniqueness}
	Let $\omega$ be a finite or infinite word and $(\alpha,\beta)\in\overline{P}$ be an alphabet. 
	If $\omega$ can be written in two ways $\gamma_1\gamma_2\dots\gamma_n$ and $\eta_1\eta_2\dots\eta_m$  where $\gamma_i,\eta_j\in\{\alpha,\beta\}$, then we have $n=m$ and $\gamma_i=\eta_i, \forall 1\leq i\leq n$, where $m,n\in\NN\cup\{\infty\}$.
\end{Prop}

\begin{proof}
	Without loss of generality we can assume that $\gamma_1=\alpha$ and $\eta_1=\beta$.
	The proposition holds trivially for the alphabet $(a,b)$.
	For any alphabet $(\alpha,\beta)\neq(a,b)\in\overline{P}$, there exist $(u,v)\in\overline{P}$ such that $(\alpha,\beta)=(uv,v) \text{ or } (u,uv)$, without loss of generality we assume it to be $(uv,v)$.
	But then by Lemma \ref{Lem3.8Har} the word $\gamma_1\gamma_2\dots=uv\dots$ will start with $v_a$ while the word $\eta_1\eta_2 \dots=v\dots$ starts with $v$, that leads to contradiction.
\end{proof}

\begin{Rmk}
	However, Proposition {\rm\ref{propuniqueness}} only proves that the representation of finite words and infinite words are unique.
	The argument for the representation of bi-infinite words follows from {\rm\cite[Lemma 3.10]{Har}} and Proposition {\rm\ref{propuniqueness}}.
\end{Rmk}

We list some properties of the renormalization $\overline{U},\overline{V}$, which are useful in the proof of the main result.

\begin{Lemma}\label{lem1}
Let $(\alpha_n,\beta_n)\in \overline{P}_n$. We have
    \begin{enumerate}
        \item $\alpha_n\beta_n=a\theta b$ with $\theta$ palindrome.\label{item1}
        \item There exists $W_n=R_1R_2\dots R_n$ with $R_i\in\{U,V\}$ such that $\alpha_n=W_n(a)$ and $\beta_n=W_n(b)$.\label{item2}
        \item Let $\omega$ be a palindrome finite word of alphabet $\{a,b\}$, then $bU(\omega)$ and $V(\omega)a$ are palindrome as well. \label{item4}
        \item $\alpha_n\alpha_n\beta_n\beta_n=a\hat{\theta}ab\hat{\theta}b$ with $\hat{\theta}$ palindrome.\label{item3}
    \end{enumerate} 
\end{Lemma}

\begin{proof}
	(\ref{item1}) and (\ref{item4}) were done in the proof of \cite[Theorem 15]{Bom}, and were restated in \cite[Remark 3.9]{Har} and \cite[Lemma 3.14]{Har}.
	(\ref{item2}) corresponds to \cite[Lemma 3.7]{Har}, which can also be regarded as a corollary of Proposition \ref{relbtrenorm}.
	(\ref{item3}) was contained in the proof of \cite[Lemma 3.15]{Har} (as a special case), but for the sake of completeness we will give the proof of (\ref{item3}):

	We use induction to prove that $W(aabb)$ is of the form $a\hat{\theta}ab\hat{\theta}b$ with $\hat{\theta}$ palindrome for any finite word $W$ of alphabet $\{U,V\}$, and then (\ref{item3}) directly follows from (\ref{item2}).
	The base case holds since we have $\hat{\theta}=\varnothing$, which is palindrome.
	Now observe that 
		\begin{align*}
			U(a\hat{\theta}ab\hat{\theta}b)&=abU(\hat{\theta})abbU(\hat{\theta})b=a\theta^{\prime}ab\theta^{\prime}b\\
			V(a\hat{\theta}ab\hat{\theta}b)&=aV(\hat{\theta})aabV(\hat{\theta})ab=a\theta^{\prime\prime}ab\theta^{\prime\prime}b
		\end{align*}
	by (\ref{item4}) we know that $\theta^{\prime}=bU(\hat{\theta}), \theta^{\prime\prime}=V(\hat{\theta})a$ are still palindrome, and this finishes the induction.
\end{proof}

In the end of this section, as was stated in the introduction, we will mainly focus on the set $m^{-1}(3)\cap l^{-1}(3)$.
The set $m^{-1}(3)$ was characterized in Proposition \ref{prop renorm}, which is the set of infinitely renormalizable bi-infinite words.
The following lemma gives a characterization of the set $m^{-1}(3)\cap l^{-1}(3)$ by renormalization.

\begin{Lemma}\label{lem4l=3}
	Let $\underline{\omega}$ be an infinitely renormalizable bi-infinite word $($so it has Markov value $3)$.
	Consider the corresponding sequence of renormalization operators $\overline{R}=\overline{R}_1\overline{R}_2\overline{R}_3\dots$, where $\overline{R}_i\in\{\overline{U},\overline{V}\}$ is the renormalization operator we applied at the $i$-th step.
	Then $\underline{\omega}$ has Lagrange value smaller than $3$ if and only if one of the operators $\overline{U},\overline{V}$ appears only finitely many times in $\overline{R}$. 
	As a consequence, $\underline{\omega}$ has Lagrange value exactly $3$ if and only if both $\overline{U}$ and $\overline{V}$ appear infinitely many times in $\overline{R}$.
\end{Lemma}

\begin{proof}
	The first remark is that we have $l(\underline{\omega})\leq m(\underline{\omega})=3$.
	
	If $\underline{\omega}$ has Lagrange value smaller than $3$, then by Markov's Theorem it is ultimately periodic to the right with period given by the Cohn tree.
	So the right side of $\underline{\omega}$ will be $\alpha_N\beta_N\alpha_N\beta_N\alpha_N\beta_N\dots$ for some $\alpha_N\beta_N\in P$ under some finite step $N$ of renormalization and in the next step the right side will be $\alpha_{N+1}\alpha_{N+1}\alpha_{N+1}\dots$ or $\beta_{N+1}\beta_{N+1}\beta_{N+1}\dots$.
	Since $\underline{\omega}$ has Markov value $3$, the only possible option is that $\underline{\omega}$ reaches degenerate case $\alpha_{N+1}^{\infty}\beta_{N+1}\alpha_{N+1}^{\infty}$ or $ \beta_{N+1}^{\infty}\alpha_{N+1}\beta_{N+1}^{\infty}$.
	So one of $\overline{U},\overline{V}$ can only appear at the first $N$ steps, which is finitely many times.
	On the other hand if one of $\overline{U},\overline{V}$ only appears for finitely many times (without loss of generality we suppose it to be $\overline{U}$), then there exist a $N_1$ big enough such that $\overline{U}$ only appears at the first $N_1$ steps.
	The same argument in Proposition \ref{prop renorm 2} shows that $\underline{\omega}$ reaches the degenerate case $\beta_{N_2}^{\infty}\alpha_{N_2}\beta_{N_2}^{\infty}$ under some finite step $N_2>N_1$, so we have $l(\underline{\omega})=l(\beta_{N_2}^{\infty}\alpha_{N_2}\beta_{N_2}^{\infty})=l(\beta_{N_2}^{\infty})<3$.
\end{proof}

\subsection{Examples of Infinitely Renormalizable Words}

Combining Proposition \ref{prop renorm 2} and Lemma \ref{lem4l=3}, we have the following characterization of the admissible bi-infinite words through renormalization, which is in the same spirit of the characterization of balanced bi-infinite words in \cite[Theorem 2.5]{Heinis}.

\begin{Prop}
	A bi-infinite word $\underline{\omega}\in\{1,2\}^{\ZZ}$ is admissible if and only if it has the following forms
	\begin{itemize}
		\item Periodic: $\underline{\omega}=\theta^{\infty}$ for some $\theta\in P$ or $\theta\in\{a,b\}$.
		\item Degenerate: $\underline{\omega}=\alpha^{\infty}\beta\alpha^{\infty}$ or $\beta^{\infty}\alpha\beta^{\infty}$ for some $(\alpha,\beta)\in\overline{P}$.
		\item Not Eventually Periodic: $\underline{\omega}$ is infinitely renormalizable but is not eventually periodic $($on neither side$)$. 
	\end{itemize}
	The periodic case corresponds to the finitely renormalizable case and the degenerate and not eventually periodic cases correspond to the infinitely renormalizable case.
	Moreover, $\underline{\omega}$ has Lagrange value $3$ if and only if it is infinitely renormalizable in the not eventually periodic case.
\end{Prop} 

Let us give some not eventually periodic examples of bi-infinite words with Markov value exactly $3$.
First we define the convergence of finite words to infinite or bi-infinite words.

\begin{Def}
	A sequence of finite words $\{w_i\}_{i\in\NN}$ converges to an infinite word $\omega$ if every prefix of $\omega$ is prefix of all but finitely many $w_i$, denoted as $\omega=\lim_{i\to\infty} w_i$.
\end{Def}

Similarly, we can define the convergence of finite words to bi-infinite words. 
Different from the infinite case, we have to set an initial position and we use a vertical bar (same notation as the cut) to denote the initial position.

\begin{Def}
	A sequence of finite words $\{w_i=w_i^1\vert w_i^2\}_{i\in\NN}$ converges to a bi-infinite word $\underline{\omega}=\omega^1\vert\omega^2=\lim_{i\to\infty}w_i^1\vert w_i^2$ if the sequence of finite words $\{w_i^2\}_{i\in\NN}$ converges to the infinite word $\omega^2$ and the sequence $\{(w_i^1)^T\}_{i\in\NN}$ converges to the infinite word $(\omega^1)^{T}$.
\end{Def}

By considering bi-infinite words as limits of some finite words, we are able to construct elements in $m^{-1}(3)\cap l^{-1}(3)$.

\begin{Prop}\label{propexample}
	Given $\overline{R}=\overline{R}_1\overline{R}_2\dots$ an infinite word of renormalization operators $\overline{U},\overline{V}$ such that both the operators $\overline{U},\overline{V}$ appears infinitely many times.
	Let $(\alpha_0,\beta_0)=(a,b)$ and $(\alpha_{n+1},\beta_{n+1})=\overline{R}_{n+1}(\alpha_{n},\beta_{n})\in\overline{P}_{n+1}$ for $n\geq 0$, then sequences $\{\beta_i\vert\alpha_i\}_{i\in\NN}$, $\{\alpha_i^T\vert\beta_i^T\}_{i\in\NN}$ will converge to some bi-infinite words with both Lagrange value and Markov value exactly $3$.
\end{Prop}

\begin{proof}
	The convergence is a consequence of the fact that $\beta_{n+1}$ ends with $\beta_n$ and $\alpha_{n+1}$ starts with $\alpha_n$, and the fact that the length of both the letters $\alpha_n,\beta_n$ tend to infinity (since both the operators $\overline{U},\overline{V}$ appear infinitely many times).
	
	We first prove that the limit word $\underline{\omega}=\lim_{i\to\infty}\beta_i\vert\alpha_i$ can be written as a non-constant word in any alphabet $(\alpha_n,\beta_n)$ for $n\in\NN$.
	This is because by Proposition \ref{propuniqueness}, words $\beta_{n^\prime}\vert\alpha_{n^\prime},\beta_{n^\prime+1}\vert\alpha_{n^\prime+1}$ can be written in a unique way in the alphabet $(\alpha_n,\beta_n)$ for any $n^\prime\geq n$.
	Moreover, $\alpha_{n^\prime}$ is a prefix of $\alpha_{n^\prime+1}$ and $\beta_{n^\prime}$ is a suffix of $\beta_{n^\prime+1}$(in alphabet $(\alpha_n,\beta_n)$).
	So the sequence $\{\beta_i\vert\alpha_i\}_{i\in\NN}$ will converge to a bi-infinite word which can be written in the alphabet $(\alpha_n,\beta_n)$.
	$\underline{\omega}$ is not a constant word in alphabet $(\alpha_n,\beta_n)$ since it contains factor $\beta_n\alpha_n$.
	Then the proposition follows directly from Proposition \ref{prop renorm 2} and Lemma \ref{lem4l=3}.
	
	As for the word $\lim_{i\to\infty}\alpha_i^T\vert\beta_i^T$, it is not hard to see that we have $\lim_{i\to\infty}\alpha_i^T\vert\beta_i^T=\underline{\omega}^T$.
	We prove that $\underline{\omega}^T$ can also be written as a non-constant word in any alphabet $(\alpha_n,\beta_n)$ for $n\in\NN$.
	First by Lemma \ref{lemtranspose} we know that $m(\underline{\omega}^T)=m(\underline{\omega})=3$.
	Then the result follows from the fact that $\underline{\omega}$ and $\underline{\omega}^T$ have the same sequence of renormalization operators.
	This is because the characteristic sequence of the word $\underline{\omega}^T$ is simply the transpose of the characteristic sequence of $\underline{\omega}$, but transpose will not change the choice of renormalization operator.
\end{proof}

\begin{Rmk}
 	Given an infinitely renormalizable bi-infinite word $\underline{\omega}=(\omega_i)_{i\in\ZZ}$ and the corresponding sequence of renormalization operators $\overline{R}=\overline{R}_1\overline{R}_2\dots$, $\overline{R}_i\in\{\overline{U},\overline{V}\}$.
 	It is natural to ask how large is the set of bi-infinite words that have precisely this sequence of renormalization operators.

 	If we define a map $p: m^{-1}(3)\rightarrow \{\overline{U},\overline{V}\}^{\NN}$ that maps every infinitely renormalizable bi-infinite word to the corresponding sequence of renormalization operators.
 	The question is that how large is the set $p^{-1}(\overline{R})$.
	It is possible to show that $p^{-1}(\overline{R})$ is uncountable if $\underline{\omega}$ is not eventually periodic, and is countable if $\underline{\omega}$ is degenerate. The uncountability follows from the well known fact that Sturmian dynamical systems are minimal \cite[Exercise 6.2.10]{substitutionsbook}, but we will give a proof for completeness.
	
	The idea is following:
	If $\underline{\omega}$ is not eventually periodic, then we consider the closure of the orbit of $\underline{\omega}$, which is $\overline{\{\sigma^n(\underline{\omega})\mid n\in\ZZ\}}$.
	We claim that the set $\overline{\{\sigma^n(\underline{\omega})\mid n\in\ZZ\}}=p^{-1}(\overline{R})$ is uncountable.
	This will follow from the fact that the set $\overline{\{\sigma^n(\underline{\omega})\mid n\in\ZZ\}}$ is a non-periodic minimal set, so in particular, it is homeomorphic to a Cantor set and contains uncountably many elements.
	
	Since the space $\{1,2\}^{\ZZ}$ is a compact metric space, it is well known that the closure of the orbit of $\underline{\omega}$ is minimal if and only if $\underline{\omega}$ is \emph{syndetically recurrent}, that is, for any $\gamma=\omega_m\omega_{m+1}\dots\omega_{m+k}\sqsubseteq\underline{\omega}$ there exists a constant $c_k\in\NN^{+}$ such that for any $n\in\ZZ$ there exist a $\tilde{n}\in\ZZ$ such that $n\leq\tilde{n}<n+c_k$ and $\gamma=\omega_{\tilde{n}}\omega_{\tilde{n}+1}\dots\omega_{\tilde{n}+k}$.
	So it suffices to show that $\underline{\omega}$ is syndetically recurrent.
	This will follow from the fact that $\underline{\omega}$ is infinitely renormalizable with both the renormalization operators $\overline{U},\overline{V}$ appearing infinitely many times $($Lemma {\rm\ref{lem4l=3}}$)$.
	For any finite word $\gamma\sqsubseteq\underline{\omega}$, we can take $N\in\NN$ big enough such that $(\alpha_N,\beta_N)=(u^{k+1}v,u^{k}v)\text{ or }(uv^k,uv^{k+1})$ with $k\geq 1$ and $\left|\gamma\right|<\min{\{\left|u\right|,\left|v\right|\}}$ for some $(u,v)\in\overline{P}$.
	It follows that the word $\gamma$ appears in any subword of size bigger than $\left|u^{k+1}v\right|+\left|uv\right|$.
	
	The containment $\overline{\{\sigma^n(\underline{\omega})\mid n\in\ZZ\}}\subseteq p^{-1}(\overline{R})$ holds because of the minimality, any large portions of $\underline{\omega}$ will appear in elements of $\overline{\{\sigma^n(\underline{\omega})\mid n\in\ZZ\}}$.
	On the other hand, the containment $p^{-1}(\overline{R})\subseteq \overline{\{\sigma^n(\underline{\omega})\mid n\in\ZZ\}}$ holds because for any $\underline{\omega}^{\prime}\in p^{-1}(\overline{R})$, any large portions of $\underline{\omega}^{\prime}$ are contained in a single letter $\alpha_n$ or $\beta_n$ for some $n$ big enough.
	It follows that $\underline{\omega}^{\prime}\in\overline{\{\sigma^n(\underline{\omega})\mid n\in\ZZ\}}$ since $\alpha_n,\beta_n$ are subwords of $\underline{\omega}$.
	
	In contrast, if $\underline{\omega}$ is degenerate, it follows from the proof of Lemma {\rm \ref{lem4l=3}} that we have $p^{-1}(\overline{R})=\{\sigma^n(\underline{\omega})\vert n\in\ZZ\}$, which is countable.
	Moreover, since $\underline{\omega}$ is eventually periodic, the set $\overline{\{\sigma^n(\underline{\omega})\mid n\in\ZZ\}}$ is countable as well.
\end{Rmk}

\section{Proof of Theorem \ref{mainthm}}\label{sec3}

Let us first explain some ideas of the proof.
Our main task is to construct uncountably many elements in $\pi(l^{-1}(3)\cap m^{-1}(3))=\{\pi(\underline{x}) \mid \underline{x}\in l^{-1}(3)\cap m^{-1}(3)\}$ with exactly $n$ bad cuts for every $n\in\NN\cup\{\infty\}$.
It follows from definition that those elements are in $Y_3(n)$. 

The idea of creating bad cuts in projections was explained right after Proposition \ref{critbadcut}, and this will be applied at the beginning of Subsection \ref{subsec3.2}.
By taking compositions of renormalizations, we first construct uncountably many elements in the sets $Y_3(\infty)$ and $Y_3(0)$. 
Then by gluing certain finite words at the beginning of elements $\omega_{X_3(0)}\in Y_3(0)$ (the finite word depends on $\omega_{X_3(0)}$), we can produce exactly $n$ bad cuts in the beginning and preserve all the good cuts in $\omega_{X_3(0)}$ (see Proposition \ref{prop4casen}).
This will prove that the set $Y_3(n)$ is uncountable for $n\in\NN^+$.

\subsection{Ideas of the construction}\label{sec3.1}

Given $x=x_1x_2\dots$ an infinite word, we can always identify $x$ with the irrational $[0;x]$, and we can naturally extend the definition of Lagrange value for irrationals to infinite words, so we will not distinguish between infinite word $x$ and irrational $[0;x]$ in the following text.

We are going to prove that there are uncountably many elements in $X_3(\infty)$ and $X_3(0)$ which are the projection of some infinitely renormalizable bi-infinite word.

We begin by considering a composition of renormalizations.
Let $(\alpha,\beta)\in\overline{P}$ be some word pair of alphabet $\{a,b\}$ and let $\widetilde{U},\widetilde{V}$ be the following actions:
	$$
	\widetilde{U}: (\alpha,\beta)\mapsto (\alpha\beta\beta,\alpha\beta\beta\beta),\quad \widetilde{V}: (\alpha,\beta)\mapsto (\alpha\beta\beta\beta,\alpha\beta\beta\beta\beta).
	$$

The renormalization by $\widetilde{U},\widetilde{V}$ can be regarded as a special case of the renormalization $\overline{U},\overline{V}$. This is because of the equalities
	\begin{equation}\label{tildeUV}
		\widetilde{U}=\overline{V}\circ\overline{U}\circ\overline{U},\quad \widetilde{V}=\overline{V}\circ\overline{U}\circ\overline{U}\circ\overline{U},
	\end{equation}
as illustrated in Figure \ref{figure3}, the red rectangle corresponds to the renormalization $\widetilde{U}$ and the blue rectangle corresponds to the renormalization $\widetilde{V}$.

\begin{figure}[h]
	\centering

\begin{tikzpicture}[
	level 1/.style={sibling distance=6cm, level distance=2cm},
	level 2/.style={sibling distance=3cm},
	edge from parent/.style={draw, -},
	every node/.style={font=\footnotesize} % Adjust font size if necessary
	]
	% Define the tree without circle shapes for nodes
	\node (A) {$(\alpha,\beta)$}
	child {node (B) {$(\alpha\beta,\beta)$}
		child {node (D) {$(\alpha\beta\beta,\beta)$}
			child{node (H) {$(\alpha\beta\beta\beta,\beta)$}
				child{node (J) {$(\alpha\beta\beta\beta\beta,\beta)$}}
				child{node (K) {$(\alpha\beta\beta\beta,\alpha\beta\beta\beta\beta)$}}
			}
			child{node (I) {$(\alpha\beta\beta,\alpha\beta\beta\beta)$}}
		}
		child {node (E) {$(\alpha\beta,\alpha\beta\beta)$}}
	}
	child {node (C) {$(\alpha,\alpha\beta)$}};
	
	% Draw rectangles around specific nodes
	\node[fit={(I)}, draw, inner sep=1pt, color=red, thick, rectangle, ] {};
	\node[fit={(K)}, draw, inner sep=1pt, color=blue, thick, rectangle] {};

    \node[right= 1mm of I] {$\tilde{U}$};
    \node[right= 1mm of K] {$\tilde{V}$};
	
	% Add labels to the edges
	\path (A) -- (B) node[midway, above] {$\overline{U}$};
	\path (A) -- (C) node[midway, above] {$\overline{V}$};
	\path (B) -- (D) node[midway, above] {$\overline{U}$};
	\path (B) -- (E) node[midway, above] {$\overline{V}$};
	\path (D) -- (H) node[midway, above] {$\overline{U}$};
	\path (D) -- (I) node[midway, above] {$\overline{V}$};
	\path (H) -- (J) node[midway, above] {$\overline{U}$};
	\path (H) -- (K) node[midway, above] {$\overline{V}$};
\end{tikzpicture}

    \caption{Renormalization $\widetilde{U}$ and $\widetilde{V}$.}
    \label{figure3}
\end{figure}

For convenience, we want all the word pairs $(\alpha,\beta)$ considered in the proof to start with $a$. Instead of starting with the root $(a,b)$, we start with $(\alpha_0,\beta_0)=(a,ab)$ and apply one of the $\widetilde{U},\widetilde{V}$ at each step of $(\alpha_n,\beta_n)$ and reach the pair $(\alpha_{n+1},\beta_{n+1})$, where  
$$
(\alpha_{n+1},\beta_{n+1})\in\{(\alpha_n\beta_n^2,\alpha_n\beta_n^3),(\alpha_n\beta_n^3,\alpha_n\beta_n^4)\}.
$$

Notice that through the above procedure, we can also draw an infinite binary tree $\widetilde{T}$: with root $(a,ab)$ and keep applying $\widetilde{U}$ to the left and $\widetilde{V}$ to the right.
Similar to previous notations, let $\widetilde{P}$ be the set of vertices of $\widetilde{T}$ and let $\widetilde{P}_n$ be the set of vertices with distance exactly $n$ to the root $(a,ab)$.
Then we have $(\alpha_n,\beta_n)\in \widetilde{P}_n$.

An important observation and also the reason why we choose the renormalization $\widetilde{U},\widetilde{V}$ is that the sequence of word pairs $\{(\alpha_n,\beta_n)\}_{n\in\NN}$ satisfies the following properties:

\begin{itemize}
	\item $\beta^{T}_{n+1}$ starts with $\beta^{T}_n$ and ends with $\alpha^{T}_n$.
	\item $\alpha_{n+1}\beta_{n+1}\beta_{n+1}$ starts with $\alpha_n\beta_n\beta_n$.
\end{itemize}
In particular, the sequences $\{\beta_n^T\}_{n\in\NN}$ and $\{\alpha_n\beta_n\beta_n\}_{n\in\NN}$ converge, so let $\omega_{X_3(0)}=\lim_{n\to\infty}\beta_n^T$ and $\omega_{X_3(\infty)}=\lim_{n\to\infty}\alpha_n\beta_n\beta_n$. 
We are going to prove that $\omega_{X_3(0)}\in X_3(0)$ and $\omega_{X_3(\infty)}\in X_3(\infty)$, that they are projections of some infinitely renormalizable bi-infinite words, and that they are uncountably many.

\subsubsection{Projection of infinitely renormalizable bi-infinite words}

Let $\pi(l^{-1}(3)\cap m^{-1}(3))=\{\pi(\underline{\omega})\mid \underline{\omega}\in l^{-1}(3)\cap m^{-1}(3)\}$ be the set of infinite words which are projections of elements in $l^{-1}(3)\cap m^{-1}(3)$.

\begin{Prop}\label{propprj}
	$\omega_{X_3(\infty)},\omega_{X_3(0)}\in \pi(l^{-1}(3)\cap m^{-1}(3))$, so in particular, $k(\omega_{X_3(0)})=k(\omega_{X_3(\infty)})=3$.
\end{Prop}

\begin{proof}
	We need to find two elements in $l^{-1}(3)\cap m^{-1}(3)$ which projects to $\omega_{X_3(0)}$ and $\omega_{X_3(\infty)}$.
	The idea is that we consider two sequences of finite words:  $\{\alpha_i^T\vert\beta_i^T\}_{i\in\NN}$ and $\{\beta_i\vert\alpha_i\beta_i\beta_i\}_{i\in\NN}$.
	
	It follows from Proposition \ref{propexample} that both of them converges, and that both the limit words $\underline{\omega}_{X_3(0)}=\lim_{i\to\infty}\alpha_i^T\vert\beta_i^T$, $\underline{\omega}_{X_3(\infty)}=\lim_{i\to\infty}\beta_i\vert\alpha_i\beta_i\beta_i$ will have Lagrange value and Markov value exactly $3$.
	This is because by (\ref{tildeUV}) the renormalization by $\widetilde{U},\widetilde{V}$ is a particular case of renormalization by $\overline{U},\overline{V}$ and both the operators $\overline{U},\overline{V}$ will appear infinitely many times.
	On the other hand, clearly we have $\pi(\underline{\omega}_{X_3(0)})=\omega_{X_3(0)}$ and $\pi(\underline{\omega}_{X_3(\infty)})=\omega_{X_3(\infty)}$, which finishes the proof.
\end{proof}

\subsubsection{Uncountably many elements}
We prove that the $\omega_{X_3(0)}$ and $\omega_{X_3(\infty)}$ we constructed are uncountably many.
Notice that we are defining two maps from infinite words of alphabet $\{\widetilde{U},\widetilde{V}\}$ to infinite words $\omega_{X_3(0)}$ or $\omega_{X_3(\infty)}$:
let $h_1: R=R_1R_2\dots\mapsto \omega_{X_3(0)}=\lim_{n\to\infty}\beta_n^T$, $h_2: R=R_1R_2\dots\mapsto \omega_{X_3(\infty)}=\lim_{n\to\infty}\alpha_n\beta_n\beta_n$, where $R_i\in\{\widetilde{U},\widetilde{V}\}$ and $(\alpha_n,\beta_n)=R_{n}R_{n-1}\dots R_1(\alpha_0,\beta_0)$ is the corresponding sequence of word pairs.

\begin{Lemma}\label{lem4uncount}
	$(\beta_n^T)^3\alpha_n^T$ is not a prefix of $(\beta_n^T)^4\alpha_n^T$, and $(\alpha_n\beta_n^2)(\alpha_n\beta_n^3)^2$ is not a prefix of $(\alpha_n\beta_n^3)(\alpha_n\beta_n^4)^2$.
\end{Lemma}

\begin{proof}
	To prove the first part, it is enough to prove that $\alpha_n^T$ is not a prefix of $\beta_n^T\alpha_n^T$. 
	Because of the Lemma \ref{Lem3.8Har}, we have
		$$
		\beta_n^T\alpha_n^T=(\alpha_n\beta_n)^T=((\beta_n)_a(\alpha_n)^b)^T=((\alpha_n)^b)^T((\beta_n)_a)^T
		$$
	so $\alpha_n^T\neq((\alpha_n)^b)^T$ is not a prefix of $\beta_n^T\alpha_n^T$.
	
	For the second part, it is enough to prove that $\alpha_n\beta_n$ is different from $\beta_n\alpha_n$. Similarly, by Lemma \ref{Lem3.8Har} we have $\alpha_n\beta_n=(\beta_n)_a(\alpha_n)^b$, which is not equal to $\beta_n\alpha_n$.
\end{proof}

\begin{Prop}\label{uncountable element}
	The maps $h_1,h_2$ are injective, so there are uncountably many infinite words $\omega_{X_3(0)},\omega_{X_3(\infty)}$.
\end{Prop}

\begin{proof}
	If $h_1$ is not injective, there exist two infinite words $R=R_1R_2\dots$ and $R^\prime=R_1^\prime R_2^\prime\dots$ of alphabet $\widetilde{U},\widetilde{V}$ such that they have the same image $\omega_{X_3(0)}$ under $h_1$.
	Let $N+1$ be the smallest index such that $R_{N+1}$ is different from $R_{N+1}^{\prime}$ and $\{(\alpha_n,\beta_n)\}_{n\in\NN},\{(\alpha_n^\prime,\beta_n^\prime)\}_{n\in\NN}$ be the corresponding sequences of word pairs. 
	We have
		$$
		\beta_{N+1}^T,(\beta_{N+1}^{\prime})^T\in\{(\beta_N^T)^3\alpha_N^T,(\beta_N^T)^4\alpha_N^T\}.
		$$
	Without loss of generality, let $\beta_{N+1}^T=(\beta_N^T)^3\alpha_N^T$ and $(\beta_{N+1}^{\prime})^T=(\beta_N^T)^4\alpha_N^T$.
	But by $\omega_{X_3(0)}=\lim_{n\to\infty}\beta_n^T=\lim_{n\to\infty}(\beta_n^{\prime})^T$ we know that $\beta_{N+1}^T$ is a prefix of $(\beta_{N+1}^{\prime})^T$, which contradicts to Lemma \ref{lem4uncount}.
	
	Similarly, by Lemma \ref{lem4uncount} we can prove $h_2$ is injective.
	So the cardinality of the elements $\omega_{X_3(0)}$ and $\omega_{X_3(\infty)}$ we found is $2^{\NN}$, which is uncountable.
\end{proof}

\subsection{Uncountably many elements in $Y_3(\infty)$}\label{subsec3.2}
We are going to prove that there are infinitely many bad cuts in $\omega_{X_3(\infty)}$.

We are using the fact that there exists an indeterminate position in the word $\alpha_n\beta_n\beta_n$, which is a prefix of the word $\omega_{X_3(\infty)}\in\pi(l^{-1}(3)\cap m^{-1}(3))$.
This is because the word $\alpha_n\alpha_n\beta_n\beta_n$ is of the form $a\hat{\theta}ab\hat{\theta}b$ with $\hat{\theta}$ palindrome.
The cut $a\hat{\theta} a\vert b\hat{\theta} b$ is a bad cut but the cut $bb\hat{\theta}a\vert  b\hat{\theta} ba$ is a good cut, so the same cut inside the word $\alpha_n\beta_n\beta_n$ is indeterminate.
As was explained in Subsection \ref{subsec2.2}, the cut is very likely to be a bad cut in $\omega_{X_3(\infty)}$.
Indeed, the following lemma shows that it is a bad cut.

\begin{Lemma}\label{lembadcut}
	There is a bad cut inside each prefix $\alpha_n\beta_n\beta_n$ of the word $\omega_{X_3(\infty)}$, more precisely, inside the first $\beta_n$ of the prefix.
\end{Lemma}

\begin{proof}
	This lemma is a consequence of the property that $\alpha_n\alpha_n\beta_n\beta_n$ is of the form $a\hat{\theta}ab\hat{\theta}b$ with $\hat{\theta}$ palindrome.
	By Proposition \ref{critbadcut} we know that the cut $a\hat{\theta}a\vert b\hat{\theta}b$ is a bad cut. 
	So if we extend $\alpha_n\alpha_n\beta_n\beta_n$ to the right to $\alpha_n\omega_{X_3(\infty)}$, the cut $a\hat{\theta}a\vert b\hat{\theta}b\dots$ inside the extended word will have value bigger than $3$.
	On the other hand, the same cut inside the word $\omega_{X_3(\infty)}$ is still a bad cut.
	This follows from the inequality
	\begin{align*}
		\lambda(\omega_{X_3(\infty)}=\dots 2\vert2 b\hat{\theta}b\dots)&=[2;2,\dots]+[0; b,\hat{\theta},b,\dots] \\
		&>[2;2,\hat{\theta},a]+[0; b,\hat{\theta},b,\dots]= \lambda(\alpha_n\omega_{X_3(\infty)}=a\hat{\theta}a\vert b\hat{\theta}b\dots)>3
	\end{align*}
	where we are using the fact that the length before the cut is odd.
	
	So there is indeed a bad cut inside the prefix $\alpha_n\beta_n\beta_n$ of the word $\omega_{X_3(\infty)}$.
	Moreover, since by the renormalization $\widetilde{U},\widetilde{V}$, we always have that the length of the word $\alpha_n$ is shorter than the length of the word $\beta_n$. 
	Consequently, the bad cut is located in the first $\beta_n$ of the $\alpha_n\beta_n\beta_n$.
\end{proof}

\begin{Prop}
	There are infinitely many bad cuts inside $\omega_{X_3(\infty)}$, so $\omega_{X_3(\infty)}\in X_3(\infty)$.
	Moreover, $\omega_{X_3(\infty)}\in Y_3(\infty)$.
\end{Prop}

\begin{proof}
	If there are only finitely many bad cuts inside $\omega_{X_3(\infty)}$, then there exist a $N\in\NN$ big enough such that all the bad cuts of $\omega_{X_3(\infty)}$ are inside the prefix $\alpha_N\beta_N\beta_N$.
	But this will contradict to the Lemma \ref{lembadcut} since there will be a bad cut inside the first $\beta_{N+1}$ of the prefix $\alpha_{N+1}\beta_{N+1}\beta_{N+1}$.
\end{proof}

\subsection{Uncountably many elements in $Y_3(0)$}

We need to prove that there are no bad cuts in $\omega_{X_3(0)}$.

\begin{Lemma}\label{inductionZ}
	Let $\alpha_n^T,\beta_n^T$ be factors of some elements in $\pi(m^{-1}(3))$ occurring at an odd position, then all the cuts inside $\alpha_n^T,\beta_n^T$ are good except the cuts at the last letter $2$.
\end{Lemma}

\begin{proof}
	First remark is that all the $\alpha_n,\beta_n$ start with letter $a$, so all the $\alpha_n^T,\beta_n^T$ end with letter $a$, which also means that end with letter $2$.
	We use induction to prove the lemma.
	In the case $n=0$, let $\alpha_0^T=a=22, \beta_0^T=ba=1122$ be factors of $\omega\in\pi(m^{-1}(3))$. 
	Notice that by Proposition \ref{oddgood}, the odd cuts are all good cuts, so we only need to check for the cut $1\vert122$, which is a good cut inside $\omega$.
	So the lemma holds for the case $n=0$.
	
	Now suppose that the lemma holds for all $n\leq r$, we are going to prove that the lemma holds for the case $n=r+1$. 
	Let $\alpha_{r+1}^T,\beta_{r+1}^T$ be factors of $\omega$, so we have 
	$$
	\alpha_{r+1}^T,\beta_{r+1}^T \in \left\{ (\beta_r^T)^2\alpha_r^T, (\beta_r^T)^3\alpha_r^T, (\beta_r^T)^4\alpha_r^T \right\}.
	$$
	By induction we know that we only need to check for the cuts at the last letter $2$ of each $\alpha_r^T$ and $\beta_r^T$ (except the last $2$ of $\alpha_{r+1}^T$ or $\beta_{r+1}^T$).
	In general, these appear only in two situations: $(\beta_r^T)^{-}a\vert\beta_r^T$ or $(\beta_r^T)^{-}a\vert\alpha_r^T$. 
	We are going to prove that those cuts are all good cuts.
	\begin{itemize}
		\item Case $1$: $(\beta_r^T)^{-}a\vert\beta_r^T$. 
		Since we start with $\beta_0=ab$, by induction we know that there exist some $(\alpha,\beta)\in\overline{P}$ such that $\beta_r=\alpha\beta$. 
		Together with (\ref{item1}) of Lemma \ref{lem1} we know that we can write $\beta_r^T$ as $b\theta_ra$, where $\theta_r$ is palindrome.
		So the cut $(\beta_r^T)^{-}a\vert\beta_r^T$ can be written as $b\theta_ra\vert b\theta_ra$, which is a good cut (Proposition \ref{critbadcut}).
		\item Case $2$: $(\beta_r^T)^{-}a\vert\alpha_r^T$.  
		Since $\beta_r^T$ always ends with $\alpha_r^T$, we can write $\beta_r^T$ as $\gamma_r^T\alpha_r^T$, where $\gamma_r$ is some finite word.
		So the cut $(\beta_r^T)^{-}a\vert\alpha_r^T$ can be written as $\gamma_r^T(\alpha_r^T)^{-}a\vert\alpha_r^T$. Same as in the case $1$, we can write $\alpha_r^T$ as $b\hat{\theta}_ra$, so the cut is of the form $\gamma_r^Tb\hat{\theta}_ra\vert b\hat{\theta}_ra$, which is also a good cut (Proposition \ref{critbadcut}).
	\end{itemize}
	This finishes the induction.
\end{proof}

\begin{Prop}\label{uncount4case0}
	There is no bad cut inside $\omega_{X_3(0)}$, so $\omega_{X_3(0)}\in X_3(0)$.
	Moreover, $\omega_{X_3(0)}\in Y_3(0)$.
\end{Prop}

\begin{proof}
	This follows from the Proposition \ref{propprj} and Lemma \ref{inductionZ}.
	Suppose that there is a bad cut inside $\omega_{X_3(0)}$.
	Since $\omega_{X_3(0)}=\lim_{n\to\infty}\beta_n^T$, the cut must be inside (and not the cut at the last letter $2$) for some $\beta_N^T$ with $N$ big enough.
	But by Lemma \ref{inductionZ}, the cut must be a good cut inside $\omega_{X_3(0)}$, which is a contradiction.
	So $\omega_{X_3(0)}\in X_3(0)$ and by Proposition \ref{propprj} it holds that $\omega_{X_3(0)}\in Y_3(0)$.
\end{proof}

\subsection{Uncountably many elements in $Y_3(n)$ for $n\in\NN_+$}

One way of proving that there are uncountably many elements in $X_3(n)$ is adding exactly $n$ coefficients equal to $3$ at the beginning of $\omega_{X_3(0)}$.

\begin{Prop}\label{adding3}
	There are exactly $n$ bad cuts inside the infinite word $\underbrace{33\dots3}_{n}\omega_{X_3(0)}$, so $\underbrace{33\dots3}_{n}\omega_{X_3(0)}\in X_3(n)$. As a consequence, the set $X_3(n)$ is uncountable.
\end{Prop}

\begin{proof}
	Clearly the cuts at the first $n$ positions (at the letter $3$) are bad cuts, so we only need to prove that there are no bad cuts inside the suffix $\omega_{X_3(0)}$.
	The proof proceeds by induction, which is similar to Lemma \ref{inductionZ}. We only need to check that every odd cut inside the suffix $\omega_{X_3(0)}$ is good, the rest part of the induction coincide with the proof of Lemma \ref{inductionZ}.
	
	Let $33\dots 3\omega_1\omega_2\dots\omega_{i-1}\vert\omega_i\omega_{i+1}\dots$ be the cut at the $n+i$-th position where $i$ is odd, the cut is good because of the inequality
	\begin{align*}
		\lambda(33\dots 3\omega_1\omega_2\dots\omega_{i-1}\vert\omega_i\omega_{i+1}\dots)&=[0;\omega_{i-1},\dots,\omega_1,3,3,\dots,3]+[\omega_i;\omega_{i+1},\dots]\\
		&<[0;\omega_{i-1},\dots,\omega_1,\omega_0,\omega_{-1},\dots]+[\omega_i;\omega_{i+1},\dots]\\
		&\leq m(\underline{\omega}_{X_3(0)})=3
	\end{align*}
	where we are using the fact that $\omega_0\in\{1,2\}$ is smaller than $3$ and the $i$-th position is an odd position.
\end{proof}

We have now completed the proof of Theorem \ref{thm1.2}.
However, it is evident that the elements we found in Proposition \ref{adding3} are not projections of elements in $m^{-1}(3)$.
To prove that the set $Y_3(n)$ is uncountable, one way is to anticipate an $\alpha_n^T$ and remove the first letter $b$ at the beginning of $\omega_{X_3(0)}$, i.e., consider the word $(\alpha_n^T)^+\omega_{X_3(0)}$.

\begin{Lemma}\label{psfix}
	For any $(\alpha_n,\beta_n)\in\widetilde{P}_n$ with $n\geq 1$, we can always write $\alpha_n=a\theta b$ and $\beta_n=a\theta^{\prime} b$ where $\theta,\theta^{\prime}$ are palindrome.
	Moreover, $\theta$ is a prefix and suffix of $\theta^{\prime}$.
\end{Lemma}

\begin{proof}
	The first part of lemma follows from (\ref{item1}) of Lemma \ref{lem1}.
	There exists $(u,v)\in\overline{P}$ such that $(\alpha_n,\beta_n)=(u,uv)$ and $u,v$ are of form $a\theta_ub,a\theta_vb$ where $\theta_u,\theta_v$ are palindrome.
	So we have $\alpha_n=a\theta_ub$ and $\beta_n=a\theta_uba\theta_vb$, the lemma follows since $\theta=\theta_u, \theta^{\prime}=\theta_uba\theta_v$ are also palindrome.
\end{proof}

\begin{Lemma}\label{n+1badcuts}
	Let $\omega\in\pi(m^{-1}(3))$ be an infinite word starting with $(\alpha_n^T)^+\beta_n^T$ or $(\beta_n^T)^+\beta_n^T$ $(n\geq1)$, then there are exactly $n+1$ bad cuts inside the prefix $(\alpha_n^T)^+$ or $(\beta_n^T)^+$ of $\omega$.
\end{Lemma}

\begin{proof}
	The proof proceeds by induction.
	For $n=1$, we have two possible cases which is $(\alpha_1,\beta_1)=(a^2bab,a^2babab)$ or $(a^2babab,a^2bababab)$.
	For simplicity we only check for the case $(\alpha_1,\beta_1)=(a^2bab,a^2babab)$: 
	we have
		$$
		\omega=(\alpha_1^T)^+\beta_1^T\dots=aba^2bababa^2\dots	\text{ or }		
		\omega=(\beta_1^T)^+\beta_1^T\dots=ababa^2bababa^2\dots
		$$
	so there are exactly $2$ bad cuts in the prefix $(\alpha_1^T)^+$ or $(\beta_1^T)^+$, which are
		$$
		a\vert ba^2bababa^2\dots,\; abaa\vert bababa^2\dots \text{ or }
	 	a\vert baba^2bababa^2\dots,\; ababaa\vert bababa^2\dots.
		$$
	
	Now suppose that the lemma holds for the case $n$ and we are going to prove for the case $n+1$.
	Let $\omega\in\pi(m^{-1}(3))$ be an infinite word starting with $(\alpha_{n+1}^T)^+\beta_{n+1}^T$ or $(\beta_{n+1}^T)^+\beta_{n+1}^T$.
	By
		$$
		(\alpha_{n+1}^T)^+,(\beta_{n+1}^T)^+ \in\left\{ (\beta_n^T)^+(\beta_n^T)^j\alpha_n^T \mid j=1,2,3 \right\}
		$$
	we know that $\omega$ also start with $(\beta_n^T)^+\beta_n^T$ and by induction we know that there are exactly $n+1$ bad cuts in the prefix $(\beta_n^T)^+$.

	We prove that the cut at the last $2$ of $(\alpha_{n+1}^T)^+$ or $(\beta_{n+1}^T)^+$ is the only bad cut in the suffix $(\beta_n^T)^j\alpha_n^T$ and this will finish the induction.
	Notice that the word $(\beta_n^T)^j\alpha_n^T$ can also be written in the form $\beta^T$ for some $\beta\in P$, so by Lemma \ref{inductionZ} we know that all the cuts inside $(\beta_n^T)^j\alpha_n^T$ are good except the cut at the last $2$.
	By Lemma \ref{psfix} we know that we can always write $\alpha_{n+1}$ as $a\theta b$ and $\beta_{n+1}$ as $a\theta^{\prime}b$, where $\theta,\theta^{\prime}$ are palindrome and $\theta$ is a prefix and suffix of $\theta^{\prime}$.
	So the cut at the last $2$ is always of the form $\hat{\theta} a\vert b \hat{\theta} \dots$, where $\hat{\theta}=\theta$ or $\theta^{\prime}$ is palindrome.
	This cut is a bad cut because of the Proposition \ref{critbadcut} and the inequality that $\lambda(\hat{\theta} a\vert b \hat{\theta} \dots)>\lambda(\dots\hat{\theta} a\vert b \hat{\theta} \dots)=3$.
	
\end{proof}

\begin{Prop}\label{prop4casen}
	For $n=1$, there is exactly $1$ bad cut inside the infinite word $a\omega_{X_3(0)}$, so $a\omega_{X_3(0)}\in X_3(1)$. For $n\geq 2$, there are exactly $n$ bad cuts inside the infinite word $(\alpha_{n-1}^T)^{+}\omega_{X_3(0)}$, so $(\alpha_{n-1}^T)^{+}\omega_{X_3(0)}\in X_3(n)$.
    
	Moreover, $a\omega_{X_3(0)}\in Y_3(1)$ and $(\alpha_{n-1}^T)^{+}\omega_{X_3(0)}\in Y_3(n)$ are projections of some infinitely renormalizable words, so the set $Y_3(n)$ with $n\geq 1$ is also uncountable.
\end{Prop}

\begin{proof}
    From the proof of Proposition \ref{propprj}, we know that $a\omega_{X_3(0)}=\pi(\sigma^{-2}(\underline{\omega}_{X_3(0)}))\in\pi(m^{-1}(3))$ since all the $\alpha_n^T$ end with $a$, and $(\alpha_{n-1}^T)^{+}\omega_{X_3(0)}=\pi(\sigma^{-k}(\underline{\omega}_{X_3(0)}))\in\pi(m^{-1}(3))$ where $k$ is the length (in the alphabet $\{1,2\}$) of the word $(\alpha_{n-1}^T)^{+}$.
	
	When $n=1$, it follows from Lemma \ref{inductionZ} that all the cuts in $a\omega_{X_3(0)}$ are good except the cut $a\vert\omega_{X_3(0)}$. The cut is a bad cut because
	$$
	\lambda(a\vert\omega_{X_3(0)})=\lambda(2\vert211R)=[0;2]+[2;1,1,R]>[2;2,1,1,R]+[0;1,1,R]=3.
	$$
	For the case $n\geq 2$, the Proposition directly follows form Lemma \ref{inductionZ} and Lemma \ref{n+1badcuts}.
	
	The set of elements of the form $a\omega_{X_3(0)}$ is uncountable, since there are uncountably many possible choices for $\omega_{X_3(0)}$. Similarly, the set of elements of the form $(\alpha_{n-1}^T)^+\omega_{X_3(0)}$ is also uncountable: even after fixing a word $\alpha_{n-1}$ (we also fix the choice of first $n-1$ steps of renormalization operator), we still have uncountably many possible choices for $\omega_{X_3(0)}$.
\end{proof}

\section{Related Questions}\label{sec4}

Remember that the definition of Lagrange value of irrationals can be written as $k(x)=\limsup_{n\to\infty}(\gamma_{n+1}+\eta_{n+1})$, and we extended the Lagrange value to bi-infinite sequences by defining $l(\underline{\omega})=\limsup_{n\to\infty}\lambda(\sigma^n(\underline{\omega}))$. If we replace the limsup by sup in the definition of Lagrange value for bi-infinite sequences, we got the Markov value for bi-infinite sequences, which is $m(\underline{\omega})=\sup_{n\in\ZZ}\lambda(\sigma^n(\underline{\omega}))$. 
A natural question is, what is the corresponding definition of “Markov value” for irrationals? More precisely, for $x\in\RR\backslash\QQ$ define
	$$
	\widetilde{m}(x)=\sup_{n\in\NN}(\gamma_{n+1}+\eta_{n+1}).
	$$
An equivalent characterization of the function $\widetilde{m}$ is
	$$
	\widetilde{m}(x)=\inf \left\{ c>0 \mid \left| x-\frac{p}{q} \right|<\frac{1}{cq^2} \text{ has no rational solution } \frac{p}{q} \right\}
	$$
which can be regarded as the best constant for badly approximable numbers. 

For the set $\widetilde{m}^{-1}(3)$, we have the containment relation
$$
X_3(0) = \widetilde{m}^{-1}(3)\cap k^{-1}(3) \subseteq X_3 \cap k^{-1}(3).
$$
Where the first equality holds because, for any $x\in k^{-1}(3)$, we have $\widetilde{m}(x)\geq 3$, and for any $x\in X_3(0)$ it holds that $\left|x-\frac{p}{q}\right|<\frac{1}{3q^2}$ has no rational solution $\frac{p}{q}$, so we have $\widetilde{m}(x)\leq 3$. 
On the other hand, if $\left|x-\frac{p}{q}\right|<\frac{1}{3q^2}$ has rational solutions, we have $\widetilde{m}(x)>3$.
So the uncountably many elements we found in $X_3(0)$ are also in $\widetilde{m}^{-1}(3)\cap k^{-1}(3)$.

Similarly, we can define a spectrum
	$$
	\widetilde{\cM}=\{\widetilde{m}(x)<\infty \mid x\in\RR\backslash\QQ\}.
	$$
It is also interesting to study the structure of $\widetilde{\cM}$.
The first observation is that, different from the Lagrange value $k$, the function $\widetilde{m}$ depends on the initial part of the continued fraction, so maybe the spectrum will not be as “good” as the Lagrange spectrum or Markov spectrum. 
The first element of $\widetilde{\cM}$ is no longer $\sqrt{5}$, in fact, it is $\widetilde{m}([1;1,1,\dots])=\frac{3+\sqrt{5}}{2}$. 
But similarly with $\cL$ and $\cM$, is the initial part of $\widetilde{\cM}$ also a discrete sequence accumulating precisely at $3$?
Does $\widetilde{\cM}$ has positive Hausdorff dimension?
Does $\widetilde{\cM}$ have Hall's ray?

\bibliographystyle{plain}
\bibliography{references}

\begin{thebibliography}{10}

\bibitem{Bom}
Enrico Bombieri.
\newblock Continued fractions and the {M}arkoff tree.
\newblock {\em Expo. Math.}, 25(3):187--213, 2007.

\bibitem{cusick1989markoff}
Thomas~W. Cusick and Mary~E. Flahive.
\newblock {\em The {M}arkoff and {L}agrange spectra}, volume~30 of {\em
  Mathematical Surveys and Monographs}.
\newblock American Mathematical Society, Providence, RI, 1989.

\bibitem{Har}
Harold Erazo, Rodolfo Gutiérrez-Romo, Carlos~Gustavo Moreira, and Sergio
  Romaña.
\newblock Fractal dimensions of the markov and lagrange spectra near $3$.
\newblock {\em Journal of the European Mathematical Society}, 2024.

\bibitem{gugu}
Nicolau Corção~Saldanha Fabio E. Brochero~Martinez, Carlos Gustavo T. de
  A.~Moreira and Eduardo Tengan.
\newblock {\em Teoria dos Números: Um passeio com primos e outros números
  familiares pelo mundo inteiro}.
\newblock Textos universitários. SBM, 6th edition, 2024.

\bibitem{substitutionsbook}
N.~Pytheas Fogg.
\newblock {\em Substitutions in dynamics, arithmetics and combinatorics},
  volume 1794 of {\em Lecture Notes in Mathematics}.
\newblock Springer-Verlag, Berlin, 2002.

\bibitem{Heinis}
Alex Heinis.
\newblock {\em Arithmetics and Combinatorics of Words of Low Complexity}.
\newblock PhD thesis, University of Leiden, 2001.

\bibitem{Hurwitz}
A.~Hurwitz.
\newblock Ueber die angen\"aherte {D}arstellung der {I}rrationalzahlen durch
  rationale {B}r\"uche.
\newblock {\em Math. Ann.}, 39(2):279--284, 1891.

\bibitem{booklagrange}
Davi Lima, Carlos Matheus, Carlos~G. Moreira, and Sergio Roma\~na.
\newblock {\em Classical and dynamical {M}arkov and {L}agrange
  spectra---dynamical, fractal and arithmetic aspects}.
\newblock World Scientific Publishing Co. Pte. Ltd., Hackensack, NJ, 2021.

\bibitem{Malyshev1981MarkovAL}
A.~V. Malyšev.
\newblock Markov and {L}agrange spectra (a survey of the literature).
\newblock {\em Zap. Nau\v cn. Sem. Leningrad. Otdel. Mat. Inst. Steklov.
  (LOMI)}, 67:5--38, 225, 1977.
\newblock Studies in number theory (LOMI), 4.

\bibitem{Markoff1879}
A.~Markoff.
\newblock Sur les formes quadratiques binaires ind\'efinies.
\newblock {\em Math. Ann.}, 15(3-4):381--406, 1879.

\bibitem{Markoff1880}
A.~Markoff.
\newblock Sur les formes quadratiques binaires ind\'efinies.
\newblock {\em Math. Ann.}, 17(3):379--399, 1880.
\newblock (S\'econd m\'emoire).

\bibitem{Moreira2016GeometricPO}
Carlos~Gustavo Moreira.
\newblock Geometric properties of the {M}arkov and {L}agrange spectra.
\newblock {\em Ann. of Math. (2)}, 188(1):145--170, 2018.

\bibitem{Reutenauer2006}
Christophe Reutenauer.
\newblock On {M}arkoff's property and {S}turmian words.
\newblock {\em Math. Ann.}, 336(1):1--12, 2006.

\end{thebibliography}

\end{document}